\documentclass[hidelinks]{article}
\newcommand{\documentdate}{3 II 2023}

\usepackage{a4wide,latexsym,graphicx,amsmath}
\usepackage{varioref}
\usepackage{xcolor}
\usepackage{ulem}
\usepackage{multirow}
\usepackage{hyperref}
\title{Multilevel Objective-Function-Free Optimization \\
  with an Application to Neural Networks Training}

\author{Serge Gratton\footnotemark[1],
        Alena Kopani\v{c}\'akov\'a\footnotemark[2],
        Philippe L. Toint\footnotemark[3]}

\newcommand{\beqn}[1]{\begin{equation}\label{#1}}
\newcommand{\eeqn}{\end{equation}}
\newcommand{\req}[1]{(\ref{#1})}
\newcommand{\ms}{\;\;\;\;}
\newcommand{\tim}[1]{\;\; \mbox{#1} \;\;}

\newcounter{algo}[section]
\renewcommand{\thealgo}{\thesection.\arabic{algo}}
\newcommand{\llem}[2]{\vspace{\baselineskip} 
\noindent\framebox[\textwidth]{\parbox{0.95\textwidth}{
\begin{lemma} \label{#1} \rm #2 \end{lemma} } } \vspace{\baselineskip} }
\newcommand{\algo}[3]{\refstepcounter{algo}
\begin{center}\begin{figure}[htbp]
\framebox[\textwidth]{
\parbox{0.95\textwidth} {\vspace{\topsep}
{\bf Algorithm \thealgo : #2}\label{#1}\\
\vspace*{-\topsep} \mbox{ }\\
{#3} \vspace{\topsep} }}
\end{figure}\end{center}}
\newcommand{\bpr}{{\bf Proof.} \hspace{1.5mm}}
\newcommand{\epr}{\hfill $\Box$ \vspace*{1em}}
\newcommand{\lthm}[2]{\vspace{\baselineskip} 
\noindent\framebox[\textwidth]{\parbox{0.95\textwidth}{
\begin{theorem} \label{#1} \rm #2 \end{theorem} } } \vspace{\baselineskip} }

\newcommand{\ii}[1]{\{ 1, \ldots, #1 \}}
\newcommand{\iiz}[1]{\{ 0, \ldots, #1 \}}
\newcommand{\iibe}[2]{\{ #1, \ldots, #2 \}}

\newcommand{\calO}{{\cal O}}

\renewcommand{\Re}{\hbox{I\hskip -2pt R}}

\newcommand{\smallRe}{\hbox{\footnotesize I\hskip -2pt R}}

\newcommand{\bigfrac}[2]{\frac{\displaystyle #1}{\displaystyle #2}}

\newcommand{\sfrac}[2]{{\scriptstyle \frac{#1}{#2}}}
\newcommand{\half}{\sfrac{1}{2}}
\newcommand{\eqdef}{\stackrel{\rm def}{=}}
\newcommand{\bigsum}{\displaystyle \sum}
\newcommand{\kap}[1]{\kappa_{\mbox{\tiny #1}}}

\newcommand{\kB}{\kap{B}}

\newcommand{\kR}{\kap{R}}
\newcommand{\flow}{f_{\rm low}}
\newcommand{\al}[1]{{\footnotesize{\sf #1}}}

\newcommand{\sgn}{{\rm sign}}
\newtheorem{theorem}{Theorem}[section]
\newtheorem{lemma}[theorem]{Lemma}

\newtheorem{corollary}[theorem]{Corollary}

\newcommand{\malg}{\al{MOFFTR}}
\newcommand{\ASTR}{\al{ASTR1}}
\newcommand{\proof}[1]{
\begin{list}{}{
\setlength{\topsep}{0.0pt}
\setlength{\partopsep}{0.0pt}
\setlength{\leftmargin}{0.025\textwidth}
\setlength{\rightmargin}{0.5\leftmargin}
\setlength{\labelwidth}{0.5\leftmargin}
\setlength{\labelsep}{0.25\leftmargin}}
\item \bpr #1 \epr \noindent
\end{list}}
\DeclareMathOperator*{\average}{average}
\DeclareMathAlphabet{\pazocal}{OMS}{zplm}{m}{n}

\newcommand{\comment}[1]{}
\newcommand{\reportonly}[1]{{#1}}

\usepackage{tikz}
\usepackage{pgfplotstable}
\usepackage{pgfplots}
\usepackage{graphicx}
\usepgfplotslibrary{groupplots}
\usetikzlibrary{matrix}
\usetikzlibrary{automata, external, shapes.geometric, fit}
\pgfplotsset{compat=1.13}

\definecolor{myblack}{RGB}{53, 53, 53}
\definecolor{myblue}{RGB}{40, 75, 99}
\definecolor{myred}{RGB}{192, 50, 33}
\definecolor{myyellow}{RGB}{255, 166, 48}
\definecolor{mywhite}{RGB}{240, 237, 238}
\definecolor{mygreen}{RGB}{0, 102, 0}

\definecolor{green1}{RGB}{9, 82, 86}
\definecolor{green2}{RGB}{8, 127, 140}
\definecolor{green3}{RGB}{6, 167, 125}
\definecolor{green4}{RGB}{79, 109, 122}
\definecolor{green5}{RGB}{192, 214, 223}
\definecolor{violet}{RGB}{26,69,131}

\definecolor{checkgreen}{rgb}{0,0.6,0}
\definecolor{phase1}{rgb}{0.008,0.655,1.000}
\definecolor{phase2}{rgb}{0.016,0.75,0.700}
\definecolor{phase3}{rgb}{0.929,0.35,0.700}
\definecolor{icsyellow}{cmyk}{0.00,0.11,0.53,0.00}

\usepackage{amsfonts}
\def \R{\hbox{I\hskip -2pt R}}

\usepackage[cal=boondoxo]{mathalpha}
\usepackage{bbm}

\newcommand{\splitatcommas}[1]{%
  \begingroup
  \begingroup\lccode`~=`, \lowercase{\endgroup
    \edef~{\mathchar\the\mathcode`, \penalty0 \noexpand\hspace{0pt plus 1em}}%
  }\mathcode`,="8000 #1%
  \endgroup
}

\date{\documentdate}

\begin{document}

\maketitle

\renewcommand{\thefootnote}{\fnsymbol{footnote}}
\footnotetext[1]{Universit\'e de Toulouse, INP, IRIT, Toulouse,
  France. Work partially supported by 3IA Artificial and Natural
  Intelligence Toulouse Institute (ANITI), French ``Investing for the
  Future - PIA3'' program under the Grant agreement ANR-19-PI3A-0004.
  Email: serge.gratton@enseeiht.fr.}
\footnotetext[2]{Division of Applied Mathematics, Brown University, Providence, USA.
  Email: alena\_kopanicakova@brown.edu.}
\footnotetext[3]{Namur Center for Complex Systems (naXys),
  University of Namur, Namur, Belgium.
  Email: philippe.toint@unamur.be.}

\begin{abstract}
A class of multi-level algorithms for unconstrained nonlinear
optimization is presented which does not require the evaluation of the
objective function. The class contains the momentum-less AdaGrad
method as a particular (single-level) instance. The choice of avoiding
the evaluation of the objective function is intended to make the
algorithms of the class less sensitive to noise, while the multi-level
feature aims at reducing their computational cost.  The evaluation
complexity of these algorithms is analyzed and their behaviour in the
presence of noise is then illustrated in the context of training deep
neural networks for supervised learning applications. 
\end{abstract}


\section{Introduction}

In many cases, optimization problems involving a large number of
variables do exhibit some kind of structure, be it sparsity of
derivatives
\cite{BeckHall18,BouaSchn98,FujiKojiNaka97,Lasse05,Toin81a,Toin81b},
specific invariance properties
\cite{ConnGoulToin92,Gay96,GoldWang93,GrieToin82a,MareRichTaka14}
or implicit spectral properties
\cite{BorzKuni06,GelmMand90,GratSartToin08,GratToin09,Korn97,LewiNash04,Nash00},
to cite three cases of interest. If the problem arises from the
discretization of an underlying infinite-dimensional setting, it has
long been known that considering different discretizations of the same
problem using different mesh sizes and carefully using them in what is
called a \textit{multi-level} or \textit{multigrid} algorithm can
bring substantial computational benefits. When this is the case, the
remarkable numerical performance is typically obtained by exploiting 
the natural hierarchy between these discretizations to successively 
eliminate the various frequency components of the error (or residual)
\cite{BrigHensMcCo00}
while, at the same time, using the fact that evaluations of functions
and derivatives are typically cheaper for coarse discretizations than
for fine ones. Multigrid methods are now a well-researched area of
numerical analysis and are viewed as a crucial tool for the solution
of linear and nonlinear systems resulting from the solution of
elliptic partial-differential equations. Similar ideas have also made
their way in nonlinear optimization, where both the MG-Opt \cite{LewiNash04,Nash00} and RMTR
\cite{GratSartToin08} multi-level frameworks have been designed to exploit the same
properties, often very successfully
(see \cite{LewiNash02,GratMoufSartToinToma10}, for instance).

Another approach of optimization for large problems has recently been
explored extensively, promoting the use of very simple
\textit{first-order methods} (see, among many others,
\cite{DuchHazaSing11,TielHint12,KingBa15,ReddKaleKuma18,ZhouChenTangYangCaoGu20,GratJeraToin22a}).
These methods have a very low
computational cost per iterations, but typically require a (sometimes
very) large number of them. Their popularity relies on several facts.
The first is that they can be shown to be convergent with a global
rate which is sometimes comparable to that of more complicated
methods.  The second is that simplicity is achieved by avoiding the
computation of the objective-function values and, most commonly, of
other derivatives than gradients (hence their name).  This in turn has
made them very robust in the presence of noise on the function and its
derivatives \cite{GratJeraToin22c}, an important feature when the problem is so large that
these quantities can only be realistically estimated (typically by
sampling) rather than calculated exactly. The context in which
optimization is performed with computing function values is sometimes
denoted by OFFO
(\underline{O}bjective-\underline{F}unction-\underline{F}ree
\underline{O}ptimization). A very large number of first-order OFFO
methods have been investigated, but, for the purpose of this paper, we
will focus on AdaGrad \cite{DuchHazaSing11}, one of the best-known
provably-convergent members of this class.  

The purpose of this paper is to demonstrate that it is\textit{
  possible, theoretically sound, and practically efficient to combine
  multi-level and OFFO algorithms}. We achieve these objectives by 
  presenting our contributions in three steps. 
\begin{itemize}
\item We first describe a novel class of multi-level OFFO algorithms
  (of which AdaGrad can be viewed as a single-level realization) 
  (Section~\ref{section:algo}).   
\item We then analyze the global rate of convergence of algorithms in
  this class, showing results matching the state of the art
  (Sections~\ref{section:convergence} and \ref{section:extensions}).   
\item We finally illustrate the use and advantages of the proposed methods in the
  context of noisy optimization problems resulting from the training
  of deep neural nets (DNNs) for supervised learning applications
  (Section~\ref{section:numerics}). 
\end{itemize}
A brief conclusion is then proposed in Section~\ref{section:conclusion}. 
\paragraph*{Notation.} 
The symbol $\|\cdot\|$ stands for the standard Euclidean
norm. If $x$ is a vector, $|x|$ is the vector whose $j$-th component
is $|x_j|$. The singular values of the matrix $M$ are denoted by $\sigma_i[M]$.

\section{The class of multilevel OFFO algorithms}\label{section:algo}

We now present an idealized multilevel OFFO framework merging ideas from
\cite{GratJeraToin22b} and \cite{GratSartToin08} and its analysis.
The problem we consider is a structured version of smooth unconstrained
optimization, that is
\beqn{prob}
\min_{x\in\smallRe^n} f(x)
\eeqn
for a twice continuously differentiable function $f$ from $\Re^n$
into $\Re$.  The problem is structured in that we assume that
we know a collection of functions which provide a ``hierarchical'' set
of approximations of the objective function $f$.  As indicated above,
this is typically the case when considering a function of a continuous
problem's discretization (the hierarchy being then given by varying the
discretization mesh) or when the objective function involves a graph
whose description may vary in its level of detail.  More specifically,
we assume that we know a collection of functions
$\{f_\ell\}_{\ell=1}^r$ such that each $f_\ell$ is a
twice-continuously differentiable function from $\Re^{n_\ell}$ to $\Re$,
the connection with our original problem being that $n_r=n$ and
$f_r(x) = f(x)$ for all $x\in \Re^n$. We also assume that, for each
$\ell=2,\ldots,r$, $f_\ell$ is ``more costly'' to evaluate/minimize
than $f_{\ell-1}$. This may be because $f_\ell$ has more variables
than $f_{\ell-1}$ (as would typically be the case if the $f_\ell$
represent increasingly finer discretizations of the same
infinite-dimensional objective), or because the structure (in terms of
partial separability, sparsity or eigenstructure) of $f_\ell$ is more
complex than that of $f_{\ell-1}$, or for any other reason. To fix
terminology, we will refer to a particular $\ell$ as a
\textit{level}. However, for $f_{\ell-1}$ to be useful at all in
minimizing $f_\ell$, there should be some relation between the
variables of these two functions. We thus assume that, for each
$\ell=2,\ldots,r$, there exist a full-rank linear operator $R_\ell$
from $\Re^{n_\ell}$ into $\Re^{n_{\ell-1}}$ (the restriction) and
another full-rank operator $P_\ell$ from $\Re^{n_{\ell-1}}$ into
$\Re^{n_\ell}$ (the prolongation) such that
\beqn{PRcond}
\omega P_\ell = R_\ell^T
\eeqn
for some known constant\footnote{For simplicity, we choose to make $\omega$
independent of $\ell$, which can always be achieved by
scaling.} $\omega >0$. In the context of multigrid algorithms, $P_\ell$ and
$R_\ell$ are interpreted as restriction and prolongation between a fine and 
a coarse grid (see
\cite{BrigHensMcCo00,Nash00,GratSartToin08,GratMoufSartToinToma10},
for instance).

Before going into further details, we 
establish an important convention on indices.  Since
we will have to identify, sometimes simultaneously, a level, an iteration
of our algorithm and a vector's component, we associate the index $\ell$ with
levels, $i$ with iterations and $j$ with components.  For instance,
$x_{\ell,i,j}$ stands for the $j$-th component of the vector $x$ at
iteration $i$ within level $\ell$. The iteration index $i$ will be
reset to zero each time a level is entered\footnote{We are well aware
that this creates some ambiguities, since a sequence of indices
$\ell,i$ can occur more than once if level $\ell$ ($\ell<r$) is used
more than once, implying the existence of more than one starting
iterate at this level. This ambiguity is resolved by the context.}.

Because our proposal is to extend the
\ASTR\ objective-function-free framework of \cite{GratJeraToin22b} to
the multilevel context, we now review the main concepts of this
algorithm and establish some notation. As the \al{TR} in the name
suggests, \ASTR\ is a trust-region optimization algorithm. This class
of algorithms is well-known to be both theoretically sound (see
\cite{ConnGoulToin00} for an in-depth presentation and \cite{Yuan15}
for a more recent survey) and practically very efficient.  
As in all trust-region methods, the next iterate at level $\ell$ is found by minimizing a
model of a level-dependent objective function $h_\ell$ within a region
where the model is trusted. In our
multilevel framework, this model can be either the (potentially
quadratic) Taylor-like
\beqn{Qmodel}
m_{\ell,i}(s) = g_{\ell,i}^Ts + \half s^TB_{\ell,i}s,
\eeqn
where $g_{\ell,i} \eqdef \nabla_x^1 h_\ell(x_{\ell,i})$ and
$B_{\ell,i}$ is a bounded Hessian approximation, or the 
lower level model defined by 
\beqn{hell-def}
h_{\ell-1}(x_{\ell-1,0}+s_{\ell-1}) 
\eqdef f_{\ell-1}(x_{\ell-1,0}+s_{\ell-1}) + v_{\ell-1}^Ts_{\ell-1},
\eeqn
where
\beqn{vell-def}
v_{\ell-1}= R_\ell g_{\ell,i} - \nabla f_{\ell-1}(x_{\ell-1,0}).
\eeqn
By convention, we set $v_r = 0$, so that, for all $s_r$,  
\beqn{fishr}
h_r(x_{r,0}+s_r) = f_r(x_{r,0}+s_r)= f(x_{r,0} + s_r)
\tim{and}
g_{r,k} = \nabla_x^1 h_r(x_{r,k}) = \nabla_x^1 f(x_{r,k}).
\eeqn
The model $h_\ell$ therefore corresponds to a modification of $f_\ell$
by a linear term that enforces the ``linear coherence'' relation 
\beqn{grel}
g_{\ell-1,0} = \nabla_x^1 h_{\ell-1}(x_{\ell-1,0}) = R_\ell g_{\ell,i}. 
\eeqn
This first-order modification \req{hell-def} is commonly used in
multigrid applications in the context of the full approximation scheme
\cite{BrigHensMcCo00}, but also in other contexts
\cite{Fish98,Nash00,AlexLewi01,LewiNash02,GratSartToin08,Kopanicakova2022e}.  We call it
``linear coherence'' because it crucially
ensures that the first-order behaviours of $h_\ell$ and $h_{\ell-1}$ are coherent in
a neighbourhood of $x_{\ell,i}$ and  $x_{\ell-1,0}$, respectively. To
see this, one checks that, if $s_\ell$ and $s_{\ell-1}$ satisfy 
$
s_\ell = P_\ell s_{\ell-1},
$
then, using \req{PRcond} and \req{grel},
\beqn{linearok}
g_{\ell,i}^Ts_\ell
= g_{\ell,i}^TP_\ell s_{i-1}
= \frac{1}{\omega} R_\ell g_{\ell,i}^Ts_{\ell-1}
= \frac{1}{\omega} g_{\ell-1,0}^Ts_{\ell-1}.
\eeqn

Once the model is defined/chosen, the typical
iteration of a trust-region method proceeds by minimizing it in a ball
centered at the current iterate, whose 
radius is adaptively computed by the algorithm, depending on past
performance. This ball can be defined in different norms, but we will focus
here on a scaled version of the ``infinity norm'' where the absolute
value of each vector component is measured individually.  In the
\ASTR\ context, the trust-region radius is computed using the size of the current
gradient and a component-wise strictly positive vector of \textit{weights},
which we do not fully define now, but which will be specified later in our
analysis. Since our multilevel algorithm is recursive,
it is also necessary to force termination at a given level when the
Euclidean norm of the prolongation of the overall step to the previous level
becomes too large, that is when the inequality
\beqn{sizeok0}
\|P_{\ell+1}(x_{\ell,i}-x_{\ell,0})\| \leq \delta_\ell
\eeqn
fails, where $\delta_{\ell}\geq 0$ 
is a bound on norm  of the step at level $\ell+1$ if $\ell < r$ or $+\infty$ otherwise.

In order to specify the algorithm, we finally define, for a vector of weights $w_\ell$ of size
$n_\ell$, the diagonal matrices
\beqn{D-def}
D(\{ w_\ell \})
\eqdef {\rm diag}\left(\frac{1}{w_{\ell,1}},\ldots,\frac{1}{w_{\ell, n_ \ell}}\right).
\eeqn
We will assume that, for  $j \in \{1, \ldots, n_\ell \}$, there exists
a constant $\varsigma_j \in (0,1]$, such that $w_{\ell, j}\geq \varsigma_j$
for each $\ell\in\ii{r}$. The \underline{M}ultilevel
\underline{O}bjective-\underline{F}unction-\underline{F}ree
\underline{T}rust-\underline{R}egion (\malg) algorithm is then
specified \vpageref{malg}. \\

\noindent Solving the original problem \req{prob} is then obtained by calling
\beqn{callmalg}
\malg(\,r,\,f,\,x_{r,0},\,\epsilon_r, \,i^{(\max)}_r, \, +\infty,\, \varsigma\, ),
\eeqn
where $i^{(\max)}_r$ is the maximum number of top level iterations
and~$\varsigma$ denotes the vector of weights' lower bounds.
In order to fix terminology, we say that iterations
at which the step is computed by Step~4 of the algorithm are \textit{Taylor
iterations}, while iterations at which the step results from the
recursive call \req{recstep} are called \textit{recursive iterations}.

\algo{malg}{\hspace*{20mm}\fbox{$x_+ = \malg(\ell,
    h_\ell,x_{\ell,0},\epsilon_\ell,i_\ell^{(\max)}, \delta_\ell, w_{\ell,0} )$}}
{
\begin{description}
\item[Step 0: Initialization. ]
  The constants $\kR\in (0,1)$, $\alpha \ge 1$,
  $\tau \in (0,1]$, $\kB  \ge 1$ and
  $\varsigma_j \in(0,1]$ ($j\in\ii{n_\ell})$ are given. Set $i=0$. 
\item[Step 1: Termination test. ]
  If $\ell<r$ and 
  \beqn{sizeok}
  \|P_{\ell+1}(x_{\ell,i}-x_{\ell,0})\| > \delta_\ell
  \eeqn
  return with $x_+=x_{\ell,i-1}$. 
  Otherwise, compute $g_{\ell,i} \eqdef \nabla_x^1 h_\ell(x_{\ell,i})$.
  If $\|g_{\ell,i}\|\leq \epsilon_\ell$ or $i=i_\ell^{(\max)}$, return with
  $x_+=x_{\ell,i}$.
\item[Step 2: Define the trust-region. ]
  Set
  \beqn{Delta-def}
  \widehat{\Delta}_{\ell,i}  = D(\{w_{\ell,i}\}) |g_{\ell,i}|
  \tim{and}
  \Delta_{\ell,i} = \left\{\begin{array}{ll}
  \widehat{\Delta}_{\ell,i} &\tim{if } \ell=r,\\
  \min\left[\bigfrac{2\delta_\ell}{\|P_{\ell+1}\|\,\|\widehat{\Delta}_{\ell,i}\|},1\right]
  \widehat{\Delta}_{\ell,i}
  &\tim{if } \ell<r.\\
  \end{array}\right.
  \eeqn
  If $i >0$, define $w_{\ell,i}\in \Re^{n_\ell}$ such that
  $w_{\ell,i,j} \geq \varsigma_j$ for $j\in\ii{n_\ell}$.
  If a Taylor step is required at iteration $i$, go to Step~4.
\item[Step 3: Recursive step. ]
  Select $w_{\ell-1,0}\in \Re^{n_{\ell-1}}$ such that $w_{\ell-1,0,j} \geq \varsigma_j$
  for $j\in\ii{n_{\ell-1}}$,
  \beqn{Rwcond}
  \fbox{\tim{``the lower-level weights are large enough''}}
  \eeqn
  \vspace*{-2mm}
  and
  \beqn{lowTR-cond}
  \left\| D(\{w_{\ell-1,0}\}) \, | R_\ell g_{\ell,i} | \right\|
  \leq \frac{\alpha \|\Delta_{\ell,i}\|}{\| P_\ell\|}.
  \eeqn   
  \vspace*{-2mm}
  If either $\ell = 1$ or
  \beqn{nodesc-cond}
     \sum_{j=1}^{n_\ell-1}\frac{ [ R_\ell g_{\ell,i}]^2_j}{w_{\ell-1,0,j}}
     < \kR \sum_{j=1}^{n_\ell}\frac{ g_{\ell,i,j}^2}{w_{\ell,i,j}}
  \eeqn
  then go to Step~4. Otherwise (i.e.\ if $\ell>1$ and \req{nodesc-cond} fails), compute
  \beqn{recstep}
  s_{\ell,i}
  =  P_\ell\Big[\malg(\ell-1, h_{\ell-1},R_\ell x_{\ell,i},
    \epsilon_{\ell-1},i_{\ell-1}^{(\max)},\alpha\|\Delta_{\ell,i}\|,w_{\ell-1,0})
    - R_\ell x_{\ell,i}\Big],
    \eeqn
  \vspace*{-2mm}
  where $h_{\ell-1}$ is given by \req{hell-def}.
\item[Step 4: Taylor step. ]
  Select a symmetric Hessian approximation $B_{\ell,i}$ such that
  \beqn{Bbound}
  \|B_{\ell,i}\| \le \kB.
  \vspace*{-2mm}
  \eeqn
  \vspace*{-2mm}
  Compute a step $s_{\ell,i}$ such that
  \beqn{sbound}
  |s_{\ell,i,j}| \le \Delta_{\ell,i,j}  \ms (j \in \ii{n_\ell}),
  \eeqn
  \vspace*{-3mm}
  and
  \beqn{GCPcond}
  g_{\ell,i}^Ts_{\ell,i} + \half s_{\ell,i}^TB_{\ell,i}s_{\ell,i}
  \le \tau\left(g_{\ell,i}^Ts_{\ell,i}^Q + \half (s_{\ell,i}^Q)^TB_{\ell,i}s_{\ell,i}^Q\right),
  \eeqn
  \vspace*{-3mm}
  where
  \beqn{sL-def}
  s_{\ell,i,j}^L = -\sgn(g_{\ell,i,j})\Delta_{\ell,i,j}   \ms (j \in \ii{n_\ell}),
  \eeqn
  \beqn{sQ-def}
  s^Q_{\ell,i} = \gamma_{\ell,i} s_{\ell,i}^L,
  \tim{ with }
  \gamma_{\ell,i} =
  \left\{ \begin{array}{ll}
  \min\left[ 1, \bigfrac{|g_{\ell,i}^Ts_{\ell,i}^L|}{(s_{\ell,i}^L)^T B_{\ell,i} s_{\ell,i}^L}\right] & \tim{if }
  (s_{\ell,i}^L)^T B_{\ell,i} s_{\ell,i}^L > 0,\\
  1 & \tim{otherwise.}
  \end{array}\right.
  \eeqn
  \vspace*{-8mm}
\item[Step 5: Update. ]
  Set
  \beqn{xupdate}
  x_{\ell,i+1} = x_{\ell,i} + s_{\ell,i},
  \eeqn
  Increment $i$ by one and return to Step~1.
\end{description}
}

Some comments are useful at this stage to further explain and motivate the details
of the algorithm.
\begin{enumerate}
\item Note that the iterations at any level are terminated in Step~1 if a
  level-dependent accuracy threshold $\epsilon_\ell$ is achieved or if
  a level-dependent maximum number of iterations $i_\ell^{(\max)}$ is
  reached.
  Also note that \req{sizeok} enforces \req{sizeok0}.
\item The componentwise trust-region radius is defined in
  \req{Delta-def}.  Observe that this choice prevents nonzero
  components of the step whenever the corresponding component of the
  gradient is zero.
  Observe also that \req{Delta-def} avoids large steps which would
  cause \req{sizeok} to fail for $i+1$. Indeed, if $\ell<r$, \req{Delta-def} and
  \req{sbound} imply that
  \[
  \|P_{\ell+1}(x_{\ell,i+1}-x_{\ell,i})\|
  = \|P_{\ell+1}s_{\ell,i}\|
  \leq 2\delta_\ell.
  \]
  Thus any $x_{\ell,i+1}$ which would violate this inequality would
  not satisfy \req{sizeok} (for $i+1$) since then
  \[
  \|P_{\ell+1}(x_{\ell,i+1}-x_{\ell,0})\|
  \ge \|P_{\ell+1}(x_{\ell,i+1}-x_{\ell,i})\| -
  \|P_{\ell+1}(x_{\ell,i}-x_{\ell,0})\|
  > 2\delta_\ell - \delta_\ell
  = \delta_\ell,
  \]
  where we used \req{sizeok} (for $i$) to derive the last inequality.
  
  The choice of model (end of Step~2) is not formally determined
  and left to the user. In a typical pattern, known in the multigrid
  literature as a ``V-cycle'', the tasks to perform at a given level
  $\ell$ is as follows. A set of standard Taylor iterations is first
  performed. Then, if $\ell> 1$ (that is the current level is not the
  lowest one) and significant progress is likely on level $\ell-1$
  (in the sense of \req{nodesc-cond} failing), one then recursively
  calls the algorithm at level $\ell-1$. A second set of Taylor
  iterations is then performed at level $\ell$.  The ``V'' shape
  suggested by the name results from the recursive application of this
  pattern at all levels.  While it is customary to specify the number
  of Taylor iterations in both sets (the ``pre-smoothing'' and
  ``post-smoothing'' in multigrid parlance), this is not required in
  \malg. Indeed the algorithm allows for a wide variety of iteration
  patterns, fixed or adaptive.
\item We next review the mechanism of Step~3 and start by noting
  that $\delta_{\ell,i}$ is the Euclidean norm of the step that would
  be allowed at iteration $(\ell,i)$, had this iteration been
  been a Taylor one (see \req{Delta-def}).
  We then select a set of weights $w_{\ell-1,0}$ to be used at the
  lower level.  Beyond being bounded below by their respective
  $\varsigma_j$, these weights have to satisfy two further conditions. The
  first, \req{Rwcond}, is expressed in a very generic way for now and
  ensures that these weights cannot be small if the weights at level
  $\ell$ are large. How this is achieved will depend on the specific
  choice of weights, as we will see below.
  The second is the seemingly obscure condition \req{lowTR-cond},
  which simply ensures that the global bound on the Euclidean norm on the total step at
  level $\ell-1$ is large enough to allow at least one iteration at
  the lower level. It states that the Euclidean length of the prolongation
  $P_\ell$ of the first lower level step  is at most some multiple
  $\alpha \geq 1$ of the Euclidean length of the (hypothetical) step
  at level $\ell$. Condition \req{nodesc-cond} then compares the
  decrease in a linear approximation of $h_{\ell-1}$ at $R_\ell
  x_{\ell,i}$ with that of the linear approximation of $h_\ell$ at
  $x_{\ell,i}$. If the former is less than a fraction $\kR$ of the
  latter, this suggests that ``significant progress at the lower level
  is unlikely'', and we then resort to continue minimization at the
  current level.  If significant progress is likely, we
  then choose to minimize $h_{\ell-1}$ at the lower level (recursive
  iteration) using the weights $w_{\ell-1,0}$, starting from $R_\ell x_{\ell,i}$ and
  within a Euclidean ball of radius $\alpha \|\Delta_{\ell,i}\|$. 

  We observe that \req{lowTR-cond} is quite easy to satisfy.  Indeed, one 
  readily checks that it is guaranteed if 
 \beqn{simple-w}
  w_{\ell-1,0,j} \ge \max\left[\varsigma_j, 
  \frac{\|P_\ell\|\,|R_\ell g_{\ell,i}|_j}{\alpha\|\Delta_{\ell,i}\|}\right].
  \eeqn
  
\item  The step at Taylor iterations is computed in Step~4 using a
  technique borrowed from the \ASTR\ algorithm \cite{GratJeraToin22b}.
  The reader familiar with trust-region theory will recognize in
  $s^Q_{\ell,i}$ a variant the ``Cauchy point'' obtained by minimizing the
  quadratic model on the intersection of the negative gradient's span
  and the trust-region (see \cite[Section~6.3.2]{ConnGoulToin00}), while
  $s_{\ell,i}^L$ is minimizer of the simpler linear model
  $g_{\ell,i}^Ts$ within the trust-region.
\item Neither the objective function $f_r$ or its approximations
  $\{f_\ell\}_{\ell=1}^{r-1}$ are ever evaluated and the optimization
  method combining them is therefore truly ``Objective-Function-Free''.
\end{enumerate}

\section{Convergence Analysis}\label{section:convergence}

Our convergence analysis is based on the following standard assumptions.

\begin{description}
\item[AS.1:] For each $\ell \in \ii{r}$, the function $f_\ell$ is
   continuously differentiable.
\item[AS.2:] For each $\ell \in \ii{r}$, the gradient $\nabla_x^1 f_\ell(x)$
   is Lipschitz continuous with Lipschitz constant $L\geq 0$, that is
   \[
   \|\nabla_x^1 f_\ell(x)-\nabla_x^1 f_\ell(y)\| \le L \|x-y\|
   \]
   for all $x,y\in \Re^n$ and all $\ell \in \ii{r}$.
\item[AS.3:] There exists a constant $\flow$ such that, for all $x$,
   $f(x)\ge \flow$.
\end{description}

\noindent
In what follows, we use the notation
\beqn{gap}
\Gamma_0 \eqdef f(x_0)-\flow
\eeqn
for the gap between the objective
function at the starting point and its lower bound.
Note that there is no assumption
that the gradients of the $f_\ell$ remain bounded.

Before considering more specific choices for the weights, we first
derive a fundamental property of the Taylor steps and strengthen
\cite[Lemma~2.1]{GratJeraToin22b} quantifying the ``linear decrease''
(that is the decrease in the simple linear model of the objective) for
Taylor iterations.

\llem{lemma:gs-tayl}{
Suppose that AS.1 and AS.2 hold. Consider a Taylor iteration $i$ at level
$\ell$. Then
\beqn{gs-tayl}
g_{\ell,i}^Ts_{\ell,i}
\le - \frac{\tau\varsigma_{\min}}{2\kB}\sum_{j=1}^{n_\ell}\frac{g_{\ell,i,j}^2}{w_{\ell,i,j}}
    + \frac{\kB}{2}\|\Delta_{\ell,i}\|^2
\eeqn
where $\varsigma_{\min} \eqdef \min_{j\in\ii{n_\ell}}\varsigma_j\in(0,1]$.
}
\proof{First note that, because of \req{sL-def} and the definition of $w_{\ell,i,j}$,
\beqn{gsL}
|g_{\ell,i}^Ts_{\ell,i}^L|
= \sum_{j=1}^{n_\ell}\frac{w_{\ell,i,j}g_{\ell,i,j}^2}{w_{\ell,i,j}^2}
\ge \sum_{j=1}^{n_\ell}\frac{\varsigma_jg_{\ell,i,j}^2}{w_{\ell,i,j}^2}
\ge \varsigma_{\min} \|s_{\ell,i}^L\|^2.
\eeqn
Suppose now that $(s_{\ell,i}^L)^TB_{\ell,i}s_{\ell,i}^L>0$ and
$\gamma_{\ell,i} < 1 $.  Then, in view of \req{GCPcond}, 
\req{sQ-def}, \req{gsL} and \req{Bbound},
\[
g_{\ell,i}^Ts_{\ell,i}^Q+\half (s_{\ell,i}^Q)^TB_{\ell,i}s_{\ell,i}^Q
= \gamma_{\ell,i} g_{\ell,i}^Ts_{\ell,i}^L+\half \gamma_{\ell,i}^2 (s_{\ell,i}^L)^TB_{\ell,i}s_{\ell,i}^L
= -\frac{(g_{\ell,i}^Ts_{\ell,i}^L)^2}{2(s_{\ell,i}^L)^TB_{\ell,i}s_{\ell,i}^L}
\le - \frac{\varsigma_{\min}|g_{\ell,i}^Ts_{\ell,i}^L|}{2\kB}.
\]
Combining this inequality with \req{sQ-def} then gives that
\beqn{mdecrease-conv}
g_{\ell,i}^Ts_{\ell,i}^Q+\half (s_{\ell,i}^Q)^TB_{\ell,i}s_{\ell,i}^Q
\le - \frac{\varsigma_{\min}}{2\kB}\sum_{j=1}^{n_\ell} \frac{g_{\ell,i,j}^2}{w_{\ell,i,j}}.
\eeqn
Alternatively, suppose that
$(s_{\ell,i}^L)^TB_{\ell,i}s_{\ell,i}^L \le 0$ or
$\gamma_{\ell,i}=1$. Then, using \req{sQ-def},
\req{sL-def} and bounds  $\kB \ge 1$ and $\varsigma_{\min}\leq 1$,
\beqn{mdecrease-nonconv}
g_{\ell,i}^Ts_{\ell,i}^Q+\half (s_{\ell,i}^Q)^TB_{\ell,i}s_{\ell,i}^Q
= g_{\ell,i}^Ts_{\ell,i}^L+\half (s_{\ell,i}^L)^TB_{\ell,i}s_{\ell,i}^L
\le \half g_{\ell,i}^Ts_{\ell,i}^L
\leq -\frac{\varsigma_{\min}}{2\kB}\sum_{j=1}^{n_\ell}\frac{g_{\ell,i,j}^2}{w_{\ell,i,j}}.
\eeqn
We thus obtain from \req{mdecrease-conv},
\req{mdecrease-nonconv} and \req{GCPcond} that
\[
g_{\ell,i,}^Ts_{\ell,i} + \half s_{\ell,i}^TB_{\ell,i}s_{\ell,i}
\le -\frac{\tau\varsigma_{\min}}{2\kB}\sum_{j=1}^{n_\ell} \frac{g_{\ell,i,j}^2}{w_{\ell,i,j}}.
\]
As a consequence, we deduce from \req{mdecrease-conv},
\req{mdecrease-nonconv}, \req{GCPcond} and \req{Bbound} that
\[
g_{\ell,i}^Ts_{\ell,i}
\le -\tau \varsigma_{\min}\sum_{j=1}^{n_\ell}\frac{g_{\ell,i,j}^2}{2 \kB
    w_{\ell,i,j}} + \half |s_{\ell,i}^TB_{\ell,i}s_{\ell,i}|
\le -\tau \varsigma_{\min} \sum_{j=1}^{n_\ell}\frac{g_{\ell,i,j}^2}{2 \kB
    w_{\ell,i,j}} + \frac{\kB}{2} \sum_{j=1}^{n_\ell} s_{\ell,i,j}^2,
\]
and \req{gs-tayl} results from \req{sbound}.
} 

\noindent
We also prove the following easy lemma.

\llem{lemma:stepsize}{Consider iteration $(\ell,i)$ in the course of
  the \malg\ algorithm.  Then
  \beqn{stepnorm}
  \|s_{\ell,i}\| \leq \alpha\|D(\{w_{\ell,i}\}) |g_{\ell,i}| \|.
  \eeqn
}

\proof{
If iteration $(\ell,i)$ is a Taylor iteration, \req{stepnorm} results
from \req{Delta-def}, \req{sbound} and the bound $\alpha \geq 1$.  If it is a recursive
iteration, we have, from \req{recstep} and  \req{sbound}, that
\[
\|s_{\ell,i}\|
= \|P_\ell(x_{\ell-1,i}-x_{\ell-1,0})\|
\leq \delta_{\ell-1}
= \alpha \|D(\{w_{\ell,i}\}) |g_{\ell,i}| \|,
\]
yielding \req{stepnorm}.
} 

\noindent
We now consider what can happen at a recursive iteration.

\llem{lemma:one-accepted}{Consider an recursive iteration $(\ell,i)$. Then
\beqn{rat-Deltas}
\| \Delta_{\ell-1,0} \| \leq \frac{\alpha}{\|P_\ell\|}\|\Delta_{\ell,i}\|,
\eeqn
the iterate $x_{\ell-1,1}$ is accepted in Step~1 of the
algorithm and at least one iteration is completed at level $\ell-1$.
}

\proof{
Because of \req{Delta-def}, \req{lowTR-cond} and the calling
sequence of \malg, we have that
\[
\| \Delta_{\ell-1,0} \|
\leq \| \widehat{\Delta}_{\ell-1,0} \|
= \| D(\{w_{\ell-1,0}\}) |g_{\ell-1,0}| \|
= \|D(\{w_{\ell-1,0}\}) |R_\ell g_{\ell,i}| \|
\leq \frac{\alpha}{\|P_\ell\|}\|\Delta_{\ell,i}\|
\]
and hence, using \req{sbound}, that 
\[
\|P_\ell(x_{\ell-1,1}-x_{\ell-1,0})\|
\leq \|P_\ell\| \,\|x_{\ell-1,1}-x_{\ell-1,0}\|
\leq \| P_\ell\| \,\| \Delta_{\ell-1,0} \| 
\leq \alpha \|\Delta_{\ell,i}\|
= \delta_{\ell-1}.
\]
Thus \req{sizeok} fails at iteration $(\ell-1,1)$, the iterate
$x_{\ell-1,1}$ is thus accepted and the desired conclusion follows.
} 

\llem{lemma:gs-bound}{Consider an recursive iteration $(\ell,i)$ and
suppose that  $i_{\ell-1} \geq 1$ iterations of the algorithm have
been completed at level $\ell-1$. Then
\beqn{e4}
\left|g_{\ell,i}^T s_{\ell,i}
- \frac{1}{\omega}
\sum_{k=0}^{i_{\ell-1}-1}g_{\ell-1,k}^Ts_{\ell-1,k}\right|
\le \frac{2i_{\ell-1}^{(\max)} L\delta_{\ell-1}^2}{\omega\sigma_{\min}[P_\ell]^2}.
\eeqn
}

\proof{
Using \req{recstep} and \req{PRcond} (see also \req{linearok}), we deduce that
\beqn{e1}
g_{\ell,i}^T s_{\ell,i}
= g_{\ell,i}^T\sum_{k=0}^{i_{\ell-1}-1} P_\ell s_{\ell-1,k}
= \frac{1}{\omega} \sum_{k=0}^{i_{\ell-1}-1} g_{\ell,i}^T R_\ell^T s_{\ell-1,k}
= \frac{1}{\omega} \sum_{k=0}^{i_{\ell-1}-1} g_{\ell-1,0}^Ts_{\ell-1,k}.
\eeqn
Now
\beqn{e2}
g_{\ell-1,0}^T s_{\ell-1,k}
= g_{\ell-1,k}^T s_{\ell-1,k} +
(g_{\ell-1,0}-g_{\ell-1,k})^Ts_{\ell-1,k}
\eqdef g_{\ell-1,k}^T s_{\ell-1,k} + \nu_{\ell-1,k},
\eeqn
where, using the Cauchy-Schwarz inequality, \req{hell-def}, AS.2 and \req{sizeok0},
\beqn{e3}
\begin{array}{lcl}
|\nu_{\ell-1,k}|
& \leq & \|g_{\ell-1,0}-g_{\ell-1,k}\|\,\|s_{\ell-1,k}\|\\*[1.6ex]
&  =   & \|\nabla_x^1f_\ell(x_{\ell-1,0})-\nabla_x^1f_\ell(x_{\ell-1,k})\|\,\|s_{\ell-1,k}\|\\*[1.6ex]
& \leq & L \| x_{\ell-1,k} - x_{\ell-1,0} \|\,\|s_{\ell-1,k}\|\\*[1.6ex]
& \leq & \bigfrac{L}{\sigma_{\min}[P_\ell]^2}\| P_\ell(x_{\ell-1,k}-x_{\ell-1,0})\|\,\|P_\ell s_{\ell-1,k}\|\\*[2ex]
& \leq & \bigfrac{2L\delta_{\ell-1}^2}{\sigma_{\min}[P_\ell]^2}.
\end{array}
\eeqn
Combining \req{e1}, \req{e2}, \req{e3} and the bound
$i_{\ell-1}\le i_{\ell-1}^{(\max)}$  then yields \req{e4}.
}

\noindent
This crucial lemma allows us to quantify what can be said of the ``linear
decrease'' at recursive iterations, and, as a consequence of
Lemma~\ref{lemma:gs-tayl}, at all iterations of the \malg\ algorithm.

\llem{lemma:gs-adag}{
Suppose that AS.1 and AS.2 hold.  
Then, for all $\ell\in\ii{r}$ and all $i\geq0$,
\beqn{gsR-adag}
g_{\ell,i}^T s_{\ell,i}
\leq - \beta_{1,r}\bigsum_{j=1}^{n}\bigfrac{g_{\ell,i,j}^2}{w_{\ell,i,j}}
     + \beta_{2,r}\bigsum_{j=1}^{n}\bigfrac{g_{\ell,i,j}^2}{w_{\ell,i,j}^2}
\eeqn
for some constants $\beta_{1,r}>0$ and $\beta_{2,r}>0$ independent of
$\ell$ and $i$.
}
\proof{
We first apply Lemma~\ref{lemma:gs-bound} to deduce that \req{e4}
holds. We also apply Lemma~\ref{lemma:one-accepted} to conclude that
\req{rat-Deltas} holds and that $i_{\ell-1}\geq 1$.
Suppose first that iteration $(\ell,i)$ is a recursive iteration 
and that $\ell$ is one plus the index of lowest level reached by the
call to \malg\ in \req{recstep}. Then each iteration of the
\malg\ algorithm at level $\ell-1$ is a Taylor iteration
and inequality \req{gs-tayl} in Lemma~\ref{lemma:gs-tayl} applies.
Thus, using \req{e4} and \req{Delta-def}, we derive that
\begin{align}
g_{\ell,i}^T s_{\ell,i}
&  =   - \bigfrac{\tau\varsigma_{\min}}{2\kB\omega}
         \bigsum_{k=0}^{i_{\ell-1}-1}\bigsum_{j=1}^{n_{\ell-1}}\bigfrac{g_{\ell-1,k,j}^2}{w_{\ell-1,k,j}}
       + \bigfrac{\kB}{2\omega}\bigsum_{k=0}^{i_{\ell-1}-1}\|\Delta_{\ell-1,k}\|^2
       + \bigfrac{2i_{\ell-1}^{(\max)}L\delta_{\ell-1}^2}{\omega\sigma_{\min}[P_\ell]^2} \nonumber\\
& \leq - \bigfrac{\tau\varsigma_{\min}}{2\kB\omega}
         \bigsum_{k=0}^{i_{\ell-1}-1}\bigsum_{j=1}^{n_{\ell-1}}\bigfrac{g_{\ell-1,k,j}^2}{w_{\ell-1,k,j}}
       + \bigfrac{\kB}{2\omega}\bigsum_{k=0}^{i_{\ell-1}-1}\frac{4\delta_{\ell-1}^2}{\sigma_{\min}[P_\ell]^2}
       + \bigfrac{2i_{\ell-1}^{(\max)}L\delta_{\ell-1}^2}{\omega\sigma_{\min}[P_\ell]^2} \nonumber\\
&   =  - \bigfrac{\tau\varsigma_{\min}}{2\kB\omega}
         \bigsum_{k=0}^{i_{\ell-1}-1}\bigsum_{j=1}^{n_{\ell-1}}\bigfrac{g_{\ell-1,k,j}^2}{w_{\ell-1,k,j}}
       + \bigfrac{2i_{\ell-1}^{(\max)}(\kB +
         L)}{\omega\sigma_{\min}[P_\ell]^2}\delta_{\ell-1}^2\label{at-rec-1}.
\end{align}
Taking now into account the fact that $i_{\ell-1}\ge 1$, ignoring now
the terms for $k\in\ii{i_{\ell-1}-1}$ in the first sum of the
right-hand side, using the definition of $\delta_\ell$ from the call
\req{recstep}, the failure of \req{nodesc-cond} and \req{Delta-def}, we obtain that 
\begin{align}
g_{\ell,i}^T s_{\ell,i}
& \leq - \bigfrac{\tau\varsigma_{\min}}{2\kB\omega}\bigsum_{j=1}^{n_{\ell-1}}\bigfrac{g_{\ell-1,0,j}^2}{w_{\ell-1,0,j}}
       + \bigfrac{2i_{\ell-1}^{(\max)}(\kB+L)}{\omega\sigma_{\min}[P_\ell]^2}\alpha^2\|\Delta_{\ell,i}\|^2\label{at-rec-2}\\
& \leq - \bigfrac{\tau\varsigma_{\min}\kR}{2\kB\omega}\bigsum_{j=1}^{n_\ell}\bigfrac{g_{\ell,i,j}^2}{w_{\ell,i,j}}
       + \bigfrac{2\alpha^2i_{\ell-1}^{(\max)}(\kB+L)}{\omega\sigma_{\min}[P_\ell]^2}\|\widehat{\Delta}_{\ell,i}\|^2\nonumber\\
& \leq - \bigfrac{\tau\varsigma_{\min}\kR}{2\kB\omega} \bigsum_{j=1}^{n_\ell}\bigfrac{g_{\ell,i,j}^2}{w_{\ell,i,j}}
       + \bigfrac{2\alpha^2i_{\ell-1}^{(\max)}(\kB+L)}{\omega\sigma_{\min}[P_\ell]^2}
         \bigsum_{j=1}^{n_\ell}\frac{g_{\ell,i,j}^2}{w_{\ell,i,j}^2}.\label{at-rec} 
\end{align}
Alternatively, if iteration $(\ell,i)$ is a Taylor iteration, \req{gs-tayl} and
\req{Delta-def} give that
\beqn{at-tayl}
g_{\ell,i}^T s_{\ell,i}
\leq - \bigfrac{\tau\varsigma_{\min}\kR}{2\kB}\bigsum_{j=1}^{n_\ell}\bigfrac{g_{\ell,i,j}^2}{w_{\ell,i,j}}
     + \kB\bigsum_{j=1}^{n_\ell}\frac{g_{\ell,i,j}^2}{w_{\ell,i,j}^2}.
\eeqn
Combining \req{at-rec} and \req{at-tayl}, we obtain that, for all
iterations at level $\ell$ (recursive and Taylor),
\begin{align*}
g_{\ell,i}^T s_{\ell,i}
\leq & ~- \frac{\kR}{\max[\omega,1]}\left[\bigfrac{\tau\varsigma_{\min}}{2\kB}\right]
\bigsum_{j=1}^{n_\ell}\bigfrac{g_{\ell,i,j}^2}{w_{\ell,i,j}} \\
& ~
+ \max\left\{\left[\kB\right],
       \bigfrac{2\alpha^2i_{\ell-1}^{(\max)}}{\max[\omega,1]\,\sigma_{\min}[P_\ell]^2}
       \Big(\left[\kB\right] +L\Big)\right\}
       \bigsum_{j=1}^{n_\ell}\frac{g_{\ell,i,j}^2}{w_{\ell,i,j}^2}.
\end{align*}
Note that the terms in square brackets correspond to the bound
\req{at-tayl} with $\ell$ replaced by $\ell-1$ (because level $\ell-1$ contains Taylor
iterations only). We may then recursively define
\beqn{beta1-def}
\beta_{1,1} \eqdef \bigfrac{\tau\varsigma_{\min}}{2\kB}
\tim{ and }
\beta_{1,\ell+1} \eqdef \frac{\kR}{\max[\omega,1]}\, \beta_{1,\ell},
\eeqn
\beqn{beta2-def}
\beta_{2,1} \eqdef \kB
\tim{ and }
\beta_{2,\ell+1} \eqdef
\max\left\{\beta_{2,\ell},
       \bigfrac{2\alpha^2i_{\ell-1}^{(\max)}}{\max[\omega,1]\,\sigma_{\min}[P_\ell]^2}
       \Big(\beta_{2,\ell} +L\Big)\right\}
\eeqn
for $\ell \in \ii{r}$ and obtain the desired conclusion.
}

\noindent
We may now deduce a central bound on the decrease of the objective
function.

\llem{lemma:decrf}{
Suppose that AS.1 and AS.2 hold.  Then, for $\ell \in \ii{r}$,
\beqn{fdecrease}
h_\ell(x_{\ell,i})-h_\ell(x_{\ell,i+1})
\ge \sum_{j=1}^{n_\ell}\frac{g_{\ell,i,j}^2}{w_{\ell,i,j}}
    \left[\beta_{1,r}-\frac{\beta_{2,r}+\half\alpha^2 L}{w_{\ell,i,j}}\right].
\eeqn
}

\proof{Successively using AS.1, AS.2, \req{stepnorm}, \req{Delta-def} and \req{gsR-adag}, we obtain that
\begin{align}
h_\ell(x_{\ell,i+1})
& \le h_ \ell(x_{\ell,i}) + g_{\ell,i}^Ts_{\ell,i}  +
      \half L\|s_{\ell,i}\|^2 \nonumber\\*[1.5ex]
& \le h_\ell(x_{\ell,i}) + g_{\ell,i}^Ts_{\ell,i}
      + \half L \alpha^2 \|D(\{w_{\ell,i}\})|g_{\ell,i}|\|^2 \nonumber\\*[1.5ex]
& = h_\ell(x_{\ell,i}) + g_{\ell,i}^Ts_{r,i}
      + \half\alpha^2 L \sum_{j=1}^{n_\ell} \frac{g_{\ell,i,j}^2}{w_{\ell,i,j}^2} \nonumber \\
& \le h_\ell(x_{\ell,i})  - \beta_{1,r}\bigsum_{j=1}^{n_\ell}\bigfrac{g_{\ell,i,j}^2}{w_{\ell,i,j}}
      + \left(\beta_{2,r} + \half\alpha^2 L \right)
      \sum_{j=1}^{n_\ell} \frac{g_{\ell,i,j}^2}{w_{\ell,i,j}^2}
\end{align}
giving \req{fdecrease}.
} 

\subsection{Divergent Weights}\label{conv-divs}

For our approach to be coherent and practical, we now have to specify
how the weights are chosen, and also make condition
\req{Rwcond} more explicit.  We start by considering ``divergent
weights'' defined as follows. We assume, in this section, that
the weights $w_{i,k}$ are chosen such that, for some power
parameter $0< \nu \leq \mu < 1$, all $i\in\ii{n}$ and some constants
$\varsigma_i\in(0,1]$,
\beqn{ws-divs}
\max[\varsigma_i,v_{i,j}]\, (i+1)^\nu \le w_{r,i,j} \le \max[\varsigma_i,v_{i,j}]\,(i+1)^\mu
\ms
(j \geq 0),
\eeqn
where, for each $i$, the $v_{i,j}$ are such that
\beqn{vikprop}
v_{i+1,j} > v_{i,j}  \tim{ implies that } v_{i+1,j} \leq |g_{r,i+1,j}|
\eeqn
and
\beqn{viklow}
v_{i,j} \geq |g_{r,i,j}| / a(i)
\eeqn
for some positive function $a(i)$ only depending on $i$. 
Using weights of the form
\beqn{divweights_maxgi}
v_{i,j} = \max_{t\in\iiz{i}}|g_{r,t,j}|
\eeqn
has resulted in good numerical performance  when applied to noisy examples
in the single-level case (see \cite{GratJeraToin22b}). 
This particular choice, referred to as the \al{MAXGI} update rule, satisfies \req{vikprop} and \req{viklow} (with $a(i)=1$). 
The associated
condition \req{Rwcond} is now specified as the requirement that
\beqn{Rwcond-divs}
\min_{j\in\ii{n_{\ell-1}}} w_{\ell-1,0,j}\geq
      \min_{j\in\ii{n_\ell}}w_{\ell,i,j}.
\eeqn
Taking \req{simple-w} into account, we see that the definition
\[
w_{\ell-1,0,j} = \max\left[\varsigma_j, \frac{\|P_\ell\|\,|R_\ell g_{\ell,i}|_j}{\alpha \| \Delta_{\ell,i} \|}, \min_{j\in\ii{n_\ell}}w_{\ell,i,j}\right]
\]
implies both \req{lowTR-cond} and \req{Rwcond-divs}.
Note that we only need \req{ws-divs}-\req{viklow} for level $r$, the necessary growth of the weights for lower levels being guaranteed by \req{Rwcond-divs}.

Lemma~\ref{lemma:gs-adag} and~\req{viklow} may be used to immediately deduce a  lower bound on the change in the objective function's value obtained at each iteration at level~$r$.

\llem{lemma:incrf}{
Suppose that AS.1 and AS.2 hold.  
Then 
\beqn{fincrease}
f(x_i)-f(x_{i+1})
= h_r(x_{r,i})-h_r(x_{r,i+1})
\ge -n( \beta_{2,r}+ \half\alpha^2 L) a(i)^2 
\eeqn

for all $i \geq 0$.
}

\proof{
Ignoring negative terms  in \req{fdecrease} with $\ell=r$ and using \req{ws-divs} and \req{viklow}, we deduce that
\begin{align*}
h_r(x_{r,i+1})
&\le h_r(x_{r,i})
    + \left(\beta_{2,r}+\half\alpha^2 L \right) \sum_{j=1}^n \frac{g_{r,i,j}^2}{w_{r,i,j}^2}\\
&\le h_r(x_{r,i})
    + \left(\beta_{2,r}+\half\alpha^2 L \right) \sum_{j=1}^n \frac{g_{r,i,j}^2}{\max[\varsigma,v_{i,j}]^2(j+1)^\nu}\\
&\le h_r(x_{r,i})
    + \left(\beta_{2,r}+\half\alpha^2 L \right) \sum_{j=1}^n \frac{g_{r,i,j}^2}{v_{i,j}^2(j+1)^\nu}
\end{align*}
\vspace{-2mm}
and \req{fincrease} follows from \req{fishr}.
} 

\noindent
We are now ready to state our main result for the \malg\ algorithm
using \req{ws-divs}-\req{viklow} and \req{Rwcond-divs}.

\lthm{theorem:divs}{
Suppose that AS.1--AS.3 hold and that the \malg\ algorithm is applied
to problem \req{prob} in a call of the form \req{callmalg}, where
the weights $w_{\ell,i,j}$ are chosen according to \req{ws-divs}-\req{viklow} and
the condition \req{Rwcond} is instantiated as \req{Rwcond-divs}.
Then, for any $\vartheta \in (0,\beta_{1,r})$, there
exists a subsequence $\{i_t\}\subseteq
\{i\}_{i_\vartheta}^\infty$ such that
\beqn{rate-divs}
\min_{k\in\iibe{i_\varsigma+1}{i_t}}\|\nabla_x^1 f(x_{r,k})\|^2
=\min_{k\in\iibe{i_\varsigma+1}{i_t}}\|g_{r,k}\|^2
\le \kappa_\diamond \frac{(i_t+1)^\mu}{i_t-i_\vartheta}
\le \frac{2\kappa_\diamond (i_\vartheta+1)}{i_t^{1-\mu}}
\eeqn
where
\beqn{istar-divs}
i_\vartheta \eqdef \left(\frac{\beta_{2,r}+\half\alpha^2
  L}{\varsigma_{\min}(\beta_{1,r}-\vartheta)}\right)^{\frac{1}{\nu}}-1,
\ms\ms
i_\varsigma \eqdef
\left(\frac{2(i_\vartheta+1)\kappa_\diamond}{\varsigma_{\min}}\right)^\frac{1}{1-\mu}
\eeqn
and
\[
\kappa_\diamond
\eqdef \frac{2}{\vartheta}\left[f(x_0)-\flow
        +n( \beta_{2,r}+ \half\alpha^2 L)\sum_{k=0}^{i_\vartheta} a(k)^2\right].
\]
}

\proof{See Appendix~\ref{proof-divs}.}

\noindent
Some comments on this result are in order.
\begin{enumerate}
\item Theorem~\ref{theorem:divs} provides useful information on the rate of convergence of the
  \malg\ algorithm beyond iteration of index $i_\varsigma$, which can
  computed \textit{a priori}. Indeed $i_\varsigma$ only depends on
  $\mu$, $\nu$ and problem's constants. If  $\{i_t\}=\{i\}_{i_\varsigma}^\infty$,
  the complexity bound to reach an iteration satisfying the accuracy
  requirement $\|g_k\| \leq \epsilon$ is then
  \beqn{order-divs}
  \calO\Big(\epsilon^{-\sfrac{2}{1-\mu}}\Big)+i_\varsigma(\mu,\nu),
  \eeqn
  which, for small values of $\nu$ and $\mu$, can be close to
  $\calO(\epsilon^{-2})$ (albeit at the price of a larger
  $i_\varsigma$). This rate is also achieved 
  by some variants of the single-level \ASTR\ algorithm (see
  \cite[Theorem~4.1]{GratJeraToin22b}.
  If $\{i_t\} \subset \{i\}_{i_\varsigma}^\infty$, one may have to
  wait for the next iteration in $\{i_t\}$ beyond \req{order-divs} for
  the gradient bound to be achieved.  Note that the rate of decay of
  the right-hand side of \req{rate-divs} depends on the index $i_t$
  (in the complete sequence) rather than on $t$ (the subsequence
  index).
  Interestingly, it is possible to prove that $\{i_t\}
  =\{i\}_{i_\vartheta}^\infty$ if one assumes that the
  objective-function's gradients remain uniformly bounded (see
  \cite[Theorem~4.1]{GratJeraToin22b}).
\reportonly{
\item As all worst-case bounds, the bound \req{rate-divs} is pessimistic. In
  this context it is especially the case because we have only considered 
  the case where all iterations before $i_\vartheta$ generate an increase 
  in the objective function which is as large as allowed by our
  assumptions. This is extremely unlikely in practice.
\item It is also possible to relax condition \req{Rwcond-divs} by only
  requiring that the right-hand side is at least a fixed fraction of
  the left-hand side. The arguments are essentially unmodifed, but
  involve yet another constant which percolates though the proofs.  We
  haven't included this possibility to avoid further notational burden.
}
\item It was proved in \cite[Theorem~4.2]{GratJeraToin22b} that the above
  complexity bound is sharp for a single-level. It is therefore also
  sharp for the multilevel case.  
\item  Note that the requirement \req{ws-divs} allows a variety of
  choices for the weights. The specific choice \req{divweights_maxgi} will be explored from
  the numerical point of view in the next section.
\end{enumerate}
  
\subsection{AdaGrad-like weights}\label{conv-adag}

Instead of focusing on \req{ws-divs}, we now consider a choice of
weights inspired by the popular AdaGrad method, where the necessary
growth in weight size is obtained by accumulating squared gradient
components. Interestingly, this will allow us to prove a complexity
result for the complete sequence of iterates (no subsequence is
involved). More specifically, given $\varsigma \in (0,1]$ and $\mu\in (0,1)$,
we define the weights for all $\ell \in \ii{r}$, all $i\geq 0$ and $j\in\ii{n_\ell}$  by
\beqn{ws-adag}
w_{\ell,i,j} \eqdef \left(\varsigma + \sum_{i=0}^j g_{\ell,i,j}^2\right)^\mu.
\eeqn
The AdaGrad weights are recovered for $\mu = \half$. The condition \req{Rwcond} is then specified as the
requirement that
\beqn{Rwcond-adag}
\|w_{\ell-1,0}\| \geq \|w_{\ell,i}\|.
\eeqn
Taking again \req{simple-w} into account, we verify that the definition
\beqn{weights_adag}
w_{\ell-1,0,j} = \max\left[1,\frac{\|w_{\ell,i}\|}{\|\widehat{w}_{\ell-1,0}\|}\right]\widehat{w}_{\ell-1,0,j}
\tim{ where }
\widehat{w}_{\ell-1,0,j} = \max\left[\varsigma_j, \frac{\|P_\ell\|\,|R_\ell g_{\ell,i}|_j}{\alpha \|\Delta_{\ell,i}\|}\right]
\eeqn
is sufficient to ensure both \req{lowTR-cond} and \req{Rwcond-adag}.

\noindent
We now state our complexity result for the variant
of the \malg\ algorithm using \req{ws-adag} and
\req{Rwcond-adag}. This result parallels
\cite[Theorem~3.2]{GratJeraToin22b} but uses the more complex
multilevel version of the linear decrease given by Lemma~\ref{lemma:decrf}.

\lthm{theorem:adag}{Suppose that AS.1--AS.3 hold and that the
\malg\ algorithm is applied to problem \req{prob} in a call of the form
\req{callmalg}, where the weights are chosen according to
\req{ws-adag} and \req{Rwcond} is instantiated as \req{Rwcond-adag}.
Then
\beqn{rate-adag}
\average_{k\in\iiz{i}}\|\nabla_x^1f(x_{r,k})\|^2
= \average_{k\in\iiz{i}}\|g_{r,k}\|^2
\le \frac{\kappa_*}{i+1},
\vspace*{-2mm}
\eeqn
where
\beqn{k3-def}
\kappa_* \eqdef
\left\{\begin{array}{ll}
\max\left[
           \varsigma,
           \left( \bigfrac{4n(\beta_{2,r}+\half L)}{\beta_{1,r}(1-2\mu)}\right)^{\sfrac{1}{\mu}},
           \bigfrac{1}{2}\left(\bigfrac{(1-2\mu)\Gamma_0}{n(\beta_{2,r}+\half L)}\right)^{\sfrac{1}{1-2\mu}}
\right]
&\tim{ if }  0<\mu<\half, \\*[2.5ex]
\max\left[
          \varsigma,\bigfrac{1}{2}e^{\frac{2\Gamma_0}{n(\beta_{2,r} + \half L)}},
          \bigfrac{\varsigma\psi^2}{2} \,\left|W_{-1}\left(-\bigfrac{1}{\psi}\right)\right|^2
\right]     
&\tim{ if }  \mu = \half, \\*[2.5ex]
\max\left[
         \varsigma,\left[\bigfrac{2^\mu}{\beta_{1,r}}
         \left(\Gamma_0+\bigfrac{n(\beta_{2,r}+L)\varsigma^{1-2\mu}}{2\mu-1}\right)\right]^{\sfrac{1}{1-\mu}}     
\right]
&\tim{ if }  \half < \mu < 1,
\end{array}\right.
\eeqn
where $\beta_{1,r}$ and $\beta_{2,r}$ are the constants defined by \req{beta1-def}
and \req{beta2-def}, respectively, 
\beqn{psi-def}
\psi\eqdef \frac{4 \max[\sfrac{3}{2}\beta_{r,1}, n(\beta_{r,2} + \half L)]}{\beta_{1,r}\sqrt{\varsigma}}
\eeqn
and  $W_{-1}$ is the second branch
of the Lambert function \cite{Corletal96}.
}

\proof{See Appendix~\ref{proof-adag}.}

\noindent
We conclude this analysis section with a few brief comments on Theorem~\ref{theorem:adag}.
\begin{enumerate}
\item It results from this theorem that, for any $\epsilon>0$, at most $\calO(\epsilon^{-2})$
  iterations of the \malg\  algorithm are needed to reduce $\|g_{r,k}\|$ below
  $\epsilon$. This result is thus stronger than that of Theorem~\ref{theorem:divs},
  unless $\{i_t\}=\{i\}_{i_\vartheta}^\infty$.  It is also equivalent (in order) to that known for single-level trust-region algorithms (see \cite[Theorem~2.3.7]{CartGoulToin22}).
\item An exact momentum-less version of AdaGrad is obtained by
  choosing $r=1$ and $\mu= \half$.  Theorem~\ref{theorem:adag}
  therefore provides a convergence analysis for both single- and multi-level versions
  of this method.
\item It was shown in \cite[Theorem~3.4]{GratJeraToin22b} that the
  global rate of convergence of the single-level ($r=1$) algorithm
  using AdaGrad-like weights cannot be better than
  $\calO(1/\sqrt{i})$.  This is therefore also the case for $r\geq 1$.
\item The bound \req{rate-adag} may be refined (although not improved
  in order) if one is ready to assume that the gradients are uniformly
  bounded.  We refer the reader to \cite{GratJeraToin22b} for a proof
  in the single-level case.
\reportonly{
\item
  As noted in this last reference, the bound involving the Lambert
  function can be replaced by a weaker but more explicit one by using the inequality
 \beqn{Lamb-bound}
 \left|W_{-1}(-e^{-x-1})\right| \leq 1+ \sqrt{2x} + x \tim{ for } w>0
 \eeqn
 \cite[Theorem~1]{Chat13}.
 Remembering that, for $\gamma_1$ and $\gamma_2$ given by \req{ggu},
 $\log\left( \frac{\gamma_2}{\gamma_1}\right) \geq \log(3) > 1$
 and setting $x = \log\left( \frac{\gamma_2}{\gamma_1}\right) -1 > 0$
 in \req{Lamb-bound} then yields that
 \[
 \left|W_{-1}\left(-\frac{\gamma_1}{\gamma_2} \right)\right|
 \le \log\left(\frac{\gamma_2}
     {\gamma_1}\right)+\sqrt{2\left(\log\left(\frac{\gamma_2}{\gamma_1}\right)-1\right)}.
 \]
 }
 \item The strict monotonicity of the weights implied by
   \req{ws-adag} can also be relaxed to provide further algorithmic
   flexibility. In turns out that \req{ws-adag} may be replaced by
   \[
   w_{\ell,i,j} \in [ \chi\,v_{\ell,i,j}, v_{\ell,i,j}]
   \tim{with}
   v_{\ell,i,j}\eqdef \left(\varsigma + \sum_{i=0}^j g_{\ell,i,j}^2\right)^\mu
   \]
for some (fixed) $\chi \in (0,1]$ without altering the nature of Theorem~\ref{theorem:adag}, in that
only the constant $\kappa_*$ is modified to explicitly involve $\chi$. Again, see
\cite{GratJeraToin22b} for a proof in the single-level case.

\end{enumerate}

\section{Extensions}\label{section:extensions}

\reportonly{
The above theory can be extended is a number of practically useful
and/or theoretically interesting ways.

\subsection{Iteration-dependent algorithmic elements}

There is much flexibility in the implementation of the \malg\ algorithm than the statement on page~\pageref{malg} suggests, in part because we have considered certain algorithm's parameters as constant.  While this is an advantage in many circumstances, it may happen that performance can be enhanced on specific problems by allowing these parameters to vary in a fixed range.  This is for instance the case for $\alpha$, the factor by which the lower-level trust-region radius can exceed the upper-level one. This is also the case of the definition of $f_\ell$ which is never used explicitly in our analysis, or of factor 2 in the right-hand side of the second part of \req{Delta-def}. Thus, it is fair to say that our analysis covers a whole class of possible algorithms.
}
\subsection{Exploiting lower-level iterations}

The multilevel theory we have presented so far is
limited\footnote{We ignore in this discussion the fact that evaluating
gradients at the lower level is typically significantly cheaper
computationally than computing them at the upper level, an
advantage sometimes crucial in practice.}
to exploiting the first iteration of each level (see the transition 
between \req{at-rec-1} and \req{at-rec-2}). This approach is fairly coarse
in the sense that it does not give any indication on why performing more 
than a single iteration at a lower level can be beneficial for convergence.  
To improve our understanding, we need to say more about how $h_{\ell-1}$ 
provides an approximation to $h_\ell$.  At iteration $(\ell,i)$, 
$\Delta_{\ell,i}$ is meant to represent the radius of the ball around 
$x_{\ell,i}$ in which the first-order Taylor model approximates 
$h_\ell$ sufficiently well. If $f_{\ell-1}$ in turn approximates 
$f_\ell$ somehow, we expect the linear decreases in $h_{\ell-1}$ 
(that is the terms $\sum_{j=1}^{n_{\ell-1}}g_{\ell-1,k,j}^2/w_{\ell-1,k,j}$) 
to be consistent with the linear decrease at level $\ell$ (that is 
$\sum_{j=1}^{n_\ell}g_{\ell-1,k,j}^2/w_{\ell-1,k,j}$) 
within the ball whose prolongation is of radius $\Delta_{\ell,i}$.  
In what follows, we consider what can be said if one assumes (or imposes) 
that condition \req{nodesc-cond} fails for all iterations $(\ell-1,k)$ rather than 
just for $(\ell-1,0)$, that is if
\beqn{new-cond} 
\sum_{j=1}^{n_\ell-1}\frac{ g_{\ell-1,k,j}^2}{w_{\ell-1,k,j}}
     \ge \kR \sum_{j=1}^{n_\ell}\frac{ g_{\ell,i,j}^2}{w_{\ell,i,j}}
     \tim{ for } k\in \iiz{i_{\ell-1}}.
\eeqn
(Compared to \req{nodesc-cond}, we may need to use a smaller $\kR$.)

Returning to Lemma~\ref{lemma:gs-adag} with this strengthened assumption, we obtain the following result.

\llem{lemma:gs-adag-2}{
Suppose that AS.1, AS.2 and \req{new-cond} hold.  
Then, for all $\ell\in\ii{r}$ and all $i\geq0$,
\beqn{gsR-adag-2}
g_{\ell,i}^T s_{\ell,i}
\leq - i_\ell^{\rm (low)}\beta_{1,r}\bigsum_{j=1}^{n}\bigfrac{g_{\ell,i,j}^2}{w_{\ell,i,j}}
     + \beta_{2,r}\bigsum_{j=1}^{n}\bigfrac{g_{\ell,i,j}^2}{w_{\ell,i,j}^2}
\eeqn
for some constants $\beta_{1,r}>0$ and $\beta_{2,r}>0$ independent of
$\ell$ and $i$, and where $i_\ell^{\rm (low)}$ is the total number of Taylor iterations from iteration $(\ell,i)$ to (and excluding) iteration $(\ell,i+1)$.
}

\proof{
We first follow the proof of Lemma~\ref{lemma:gs-adag} up to
\req{at-rec-1}. Then, instead of ignoring terms to obtain
\req{at-rec-2}, we keep them and use \req{new-cond} and \req{Delta-def} to deduce that
\begin{align}
g_{\ell,i}^T s_{\ell,i}
& \leq - i_{\ell-1}\bigfrac{\tau\varsigma_{\min}}{2\kB\omega}\bigsum_{j=1}^{n_{\ell-1}}\bigfrac{g_{\ell-1,0,j}^2}{w_{\ell-1,0,j}}
       + \bigfrac{2i_{\ell-1}^{(\max)}(\kB+L)}{\omega\sigma_{\min}[P_\ell]^2}\alpha^2\|\Delta_{\ell,i}\|^2\nonumber\\
& \leq - i_{\ell-1} \bigfrac{\tau\varsigma_{\min}\kR}{2\kB\omega}\bigsum_{j=1}^{n_\ell}\bigfrac{g_{\ell,i,j}^2}{w_{\ell,i,j}}
       + \bigfrac{2\alpha^2i_{\ell-1}^{(\max)}(\kB+L)}{\omega\sigma_{\min}[P_\ell]^2}\|\widehat{\Delta}_{\ell,i}\|^2\nonumber\\
& \leq - i_{\ell-1}  \bigfrac{\tau\varsigma_{\min}\kR}{2\kB\omega} \bigsum_{j=1}^{n_\ell}\bigfrac{g_{\ell,i,j}^2}{w_{\ell,i,j}}
       + \bigfrac{2\alpha^2i_{\ell-1}^{(\max)}(\kB+L)}{\omega\sigma_{\min}[P_\ell]^2}
         \bigsum_{j=1}^{n_\ell}\frac{g_{\ell,i,j}^2}{w_{\ell,i,j}^2},\label{at-rec-b}
\end{align}
where, as in Lemma~\ref{lemma:gs-adag}, we used the failure of
\req{nodesc-cond} and \req{Delta-def} to obtain the last two
inequalities.
As in Lemma~\ref{lemma:gs-adag}, \req{at-tayl} still holds if iteration $(\ell,i)$ 
is a Taylor iteration (in which case $i_{\ell-1}=0$).
Combining \req{at-rec-b} and \req{at-tayl}, we obtain that, for all
iterations at level $\ell$,
\begin{align*}
g_{\ell,i}^T s_{\ell,i}
\leq & ~-\max[1,i_{\ell-1}] \frac{\kR}{\max[\omega,1]}\left[\bigfrac{\tau\varsigma_{\min}}{2\kB}\right]
       \bigsum_{j=1}^{n_\ell}\bigfrac{g_{\ell,i,j}^2}{w_{\ell,i,j}} \\
       & ~
       + \max\left\{\left[\kB\right],
       \bigfrac{2\alpha^2i_{\ell-1}^{(\max)}}{\max[\omega,1]\,\sigma_{\min}[P_\ell]^2}
       \Big(\left[\kB\right] +L\Big)\right\}
       \bigsum_{j=1}^{n_\ell}\frac{g_{\ell,i,j}^2}{w_{\ell,i,j}^2}.
\end{align*}
We may then recursively  define $\beta_{1,\ell}$ and $\beta_{1,\ell}$
using \req{beta1-def} and \req{beta2-def}
and \req{gsR-adag-2} finally follows from the definition of $i_\ell^{\rm (low)}$.
}

\noindent
Observe that \req{gsR-adag-2} is the same as \req{gsR-adag}, except
that $\beta_{1,r}$ is now multiplied by $i_\ell^{\rm (low)}$. This
modification percolates through all proofs, resulting in improved\footnote{In 
particular, by offsetting the effect of the $i_{\ell-1}^{(\max)}$ constants in $\beta_{2,r}$.}
constants in \req{istar-divs} and \req{k3-def}. 
We finally note that requiring \req{new-cond} can be viewed as one way to improve the 
balance between negative and positive terms in \req{gsR-adag-2}, but may not be the only one.

\subsection{Weak coherence}\label{relaxed}

  It is possible to relax somewhat the linear coherence requirement 
  between high and low levels models. Examination of the above theory shows
  that it is only used in \req{e1}.  If we were to assume that
  $
  \omega P_\ell^T= R_\ell + E_\ell
  $
  instead of \req{PRcond}, then it is easy to verify that an error matrix $E_\ell$
  satisfying
  \beqn{weak-coherence}
  \|E_\ell g_{\ell,i} \| \leq \kappa_E \delta_{\ell-1}
  \eeqn
  for some fixed $\kappa_E \ge 0$ ensures that, for each $t$ in \req{e1},
  \[
  g_{\ell,i}^T P_\ell s_{\ell-1,t}
  = \frac{1}{\omega} (R_\ell g_{\ell,i})^T s_{\ell-1,t} + \frac{\kappa_E\delta_{\ell-1}\|s_{\ell-1,t}\|}{\omega}
  \le \frac{1}{\omega} g_{\ell-1,i} ^T s_{\ell-1,t}+ \frac{2\kappa_E}{\omega\sigma_{\min}[P_\ell]}\,\delta_{\ell-1}^2,
  \]
  where we used \req{sizeok0} to derive the last inequality. This adds a term 
  in $\calO(\delta_{\ell-1}^2)$ in the right-hand side of \req{e4}, allowing 
  the argument to be continued with different constants.  The condition
  \req{weak-coherence} is implementable because $\delta_{\ell-1}$ is known 
  before $R_\ell$ or $P_\ell$ is used at iteration $(\ell,i)$ and may be quite 
  lax in the early iterations where $\|g_{\ell,i}\|$ is still
  relatively large.
That linear coherence often only needs to be preserved approximately is of 
particular relevance when the algorithm is applied to problems whose 
gradient is noisy.  In that case, insisting on exact linear coherence would merely 
propagate the error in the gradient at the upper level to the lower level, 
which is clearly undesirable.

\section{Numerical results}\label{section:numerics}

In this section, we illustrate the numerical performance of the
proposed \malg\ algorithms in the context of deep neural networks'
training with a particular focus on supervised learning
applications. Let~${\pazocal{D} = \{ (y_s, c_s) \}_{s=1}^{n_s}}$ be a
dataset of labeled data, where~$y_s \in \R^{n_{in}}$ represents input
features and~${c_s \in \R^{n_{out}}}$ denotes a desirable target. Our
goal is to learn the parameters of DNNs, such that they can
approximate~$c_s$ for a given~$y_s$.  
In what follows, we exploit a continuous-in-depth approach to DNNs,
the forward propagation of which can be interpreted as a
discretization of a nonlinear ordinary differential equation
(ODE)~\cite{queiruga2020continuous,chang2017multi,weinan2017proposal}.  
This approach allows us to construct a multilevel hierarchy and
transfer operators required by the \malg\ framework in a fairly
natural way (see Section~\ref{sec:algo}). 

Using a continuous-in-depth approach, the supervised learning problem
can be formulated as the following continuous optimal control
problem~\cite{haber2017stable}: 
\begin{align}
 &\underset{Q, q, \theta, W_T, b_T}{\text{min}}
 \ \ \ \frac{1}{n_{s}} \sum_{s=1}^{n_{s}} \mathcal{h}(\pazocal{P}(W_T {q}_s(T) + b_T), c_s) 
+  \int\limits_{0}^{T} \pazocal{R}({\theta}(t) ) \ dt 
 +  \pazocal{S}({W}_T, b_T) ,
 \nonumber \\  
 & \text{subject to} \quad \partial_t {q}_s(t)
    = \pazocal{F}({q}_s(t), {\theta}(t)), \qquad \forall t \in (0,T), \label{eq:cts_problem} \\
& \quad \quad \quad \quad \quad \quad {q}_s(0) = Qy_s, \nonumber
\end{align}
where $q$ denotes time-dependent states from $\R$ into $\R^{n_{fp}}$
and ${\theta}$ denotes the time-dependent control parameters from
$\R$ into  $\R^{n_{c}}$. The constraint in~\eqref{eq:cts_problem}
continuously transforms an input feature~${y_s}$ into final
state~${q}_s(T)$, defined at the time $T$. This is achieved in two steps. 
Firstly, the input $y_s$ is mapped into the dimension of the dynamical
system as~${q}_s(0) = Qy_s$, where~${Q} \in \R^{{n_{fp}} \times n_{in}}$. 
Secondly, the nonlinear transformation of the features is performed
using a "residual block"~$\pazocal{F}$ from $\R^{n_{fp}} \times \R^{n_{c}}$ into $\R^{n_{fp}}$.
We consider two types of such blocks: dense and convolutional. 
A dense residual block is defined
as~$\pazocal{F}({q}_s(t), {\theta}(t)):=\sigma(W(t) q_s(t) + b(t))$,
where ${\theta(t) = (\text{flat}(W(t)), \text{flat}(b(t)))}$,
$\sigma$ is an activation function from~$\R^{{n_{fp}}}$
into~$\R^{{n_{fp}}}$, ${b(t) \in \R^{n_{fp}}}$ is the "bias"
and ${W(t) \in \R^{n_{fp} \times n_{fp}}}$ is a dense matrix. 
A convolutional residual block has the form
${\pazocal{F}({q}_s(t), {\theta}(t)):= \sigma(\text{BN}(t, W(t) q_s(t)  + b(t)))}$,
where~$W(t)$ now stands for a sparse convolutional operator and
$\text{BN}$ denotes continuous-in-time
batch-normalization~\cite{ioffe2015batch, queiruga2020continuous} from
$\R \times \R^{n_{fp}}$ into $\R^{n_{fp}}$.

The objective function in~\eqref{eq:cts_problem} is defined such that
the deviation, measured by the loss function~$\mathcal{h}$ from
$\R^{n_{out} \times n_{out}}$ into $\R$,  between the desirable
target~$c_s$ and predicted output~$\hat{c}_s \in \R^{n_{out}}$ is
minimized.  Here, the predicted output is obtained as~$\hat{c}_s :=
\pazocal{P}(W_T {q}_s(T) + b_T) \in \R^{n_{out}}$, where~$\pazocal{P}$
denotes a hypothesis function from~$\R^{n_{out}} $
into~$\R^{n_{out}}$. The linear operators $W_T \in \R^{n_{out} \times
  n_{fp}}$, and $b_T \in \R^{n_{fp}}$ are used to perform an affine
transformation of the extracted features~${q}_s(T)$, i.e.,~features
obtained as an output of the dynamical system at time $T$.  The
regularizers~$\pazocal{R}$ and $\pazocal{S}$ with parameters~$\beta_1,
\beta_2 > 0$ are defined as follows. A Tikhonov regularization is used
to penalize the magnitude of $W_T$ and $b_T$, i.e.,~$\pazocal{S}(W_T,
b_T):= \frac{\beta_1}{2} \| W_T \|^2_F + \frac{\beta_1}{2} \| b_T
\|^2$, where $\| \cdot \|_F$ denotes the Frobenius norm.  For the
time-dependent control parameters, we use~$\pazocal{R}(\theta(t)) :=
\frac{\beta_1}{2} \| \theta(t) \|^2  + \frac{\beta_2}{2} \| \partial_t
\theta(t) \|^2$, where the second term ensures that the parameters
vary smoothly in time, see~\cite{haber2017stable} for details.   

To solve the problem~\eqref{eq:cts_problem} numerically, we follow the
first-discretize-then-optimize approach.  The discretization is
performed using equidistant grid ${0 = \tau_0 < \cdots < \tau_{K-1} =  T}$,
consisting of $K$ points. The states and controls are then
approximated at a given time~$\tau_k$ as~${{q}_k \approx {q}(\tau_k)}$,
and ${{\theta}_k \approx {\theta}(\tau_k)}$, respectively.
Note, each~$\theta_k$ and~${q}_k$ now corresponds to parameters
and states associated with the $k$-th layer of the DNN.
Our time-discretization uses the forward Euler scheme, which gives
rise to the well-known ResNet architecture with identity
skip connections~\cite{he2016identity}.  Alternatively, one could
employ more advanced, and perhaps more numerically stable, time
integration schemes, see for example~\cite {haber2017stable}.
In the case of the explicit Euler scheme considered here, we ensure
numerical stability by employing a sufficiently small time-step~$\Delta_t = T/(K-1)$. 

\subsection{Implementation and algorithmic setup}
\label{sec:algo}

Our implementation of ResNets is based on the deep-learning library
Keras~\cite{chollet2015keras}, while the \malg\ framework is implemented
using the library NumPy. We consider four different variants of
first-order \malg\ algorithms (i.e.\ $B_{\ell,i} = 0$ for all
$(\ell,i)$). The first variant employs divergent weights, specified by
the \al{MAXGI} update rule given by \req{ws-divs}
and\eqref{divweights_maxgi} with ${\mu=0.1}$.  All other variants use
AdaGrad-like weights, specified by~\eqref{ws-adag} and
\req{weights_adag}, with ${\mu \in \{0.1, 0.5, 0.9\}}$.  The selected
update rules are used to update weights at all levels.
The~\malg\ algorithms are implemented as a V-cycle with one
pre-smoothing step and zero post-smoothing steps.  For the ResNets
with dense residual blocks, we perform $10$ iterations on the lowest
level, i.e.,~$i_1^{(\max)}=10$ and employ~$\kappa_R=0.01$, and
$\alpha=5$.  For the convolutional ResNets, we use~$i_1^{(\max)}=5$ ,
$\kappa_R=0.001$ and~$\alpha=25$.  Parameters~$\nu$ and~$\varsigma$
are set as~${\nu=0.1}$, and ${\varsigma=0.01}$ for all numerical
examples. Moreover, we take our discussion of Section~\ref{relaxed}
into account and do not impose the first-order coherence
relation~\eqref{grel}, as we apply the \malg\ framework in stochastic
settings where subsampling noise is present.  In our realization of
the \malg \ algorithm, the subsampled derivatives are used at
all~$(\ell, i)$. 

The hierarchy of objective functions~$\{f\}_{\ell=1}^r$ required by
the \malg\ framework is obtained by discretizing the
problem~\eqref{eq:cts_problem} with varying discretization
parameter~$\Delta_t$.  Each~$f_{\ell}$ is then associated with a
network of different depth.  Unless stated otherwise, all numerical
examples considered below take advantage of three levels, which we
obtain using uniform refinement with a factor of two.  The
operators~$\{ P \}_{\ell=1}^{r-1}$ are constructed using piecewise
linear interpolation in 1D (in time), see~\cite{Kopanicakova_2020c,
  kopanicakova2022globally} for details.  Note that similar approaches
for assembly of prolongation operators were also employed in the
context of multilevel parameter initialization
in~\cite{haber2018learning, cyr2019multilevel, chang2017multi}.  We
define the restriction operators $\{ R \}_{\ell=1}^{r-1}$ from
\req{PRcond}, choosing~$\omega=\half$ and~$\omega=1$ for networks with
dense and convolutional residual blocks, respectively.

In order to assess the performance of the \malg \ method, we provide a
comparison with the single-level \ASTR\ methods, which we obtain by
calling the corresponding \malg\ algorithm with~$r=1$.  Our comparison
also includes the baseline stochastic gradient
(\al{SGD})~\cite{RobbMunr51}, \al{ADAM}~\cite{KingBa15} and
\al{AdaGrad}\footnote{AdaGrad is obtained by calling the
\malg\ algorithm with $r=1$, weights given by \eqref{ws-adag} and $\mu
= \half$.} methods.  The learning rate of all methods is chosen by
thorough hyper-parameter search, performed individually for each
dataset, network, and batch size.  More precisely, we consider
learning rates from the set $\splitatcommas{ \{ 0.0001, 0.00025,
  0.0005, 0.00075, 0.001, 0.0025, 0.005, 0.0075, 0.01, 0.025, 0.05,
  0.075, 0.1, 0.25, 0.5, 0.75, 1.0\}}$.  The learning rate, which gave
rise to the best generalization results (averaged over five
independent runs), is then used in the presented numerical
experiments.  In the context of the \malg\ method, the same learning
rate is employed on all levels.

In what follows, the comparison between single and multilevel methods
is performed by analyzing their dominant computational cost, i.e.,
that associated with gradient evaluations.  Let $C_r$ be a
computational cost associated with an evaluation of the gradient on
the uppermost level using a full dataset~$\pazocal{D}$.  Using the
definition of~$C_r$ and taking advantage of the fact that the cost of
the back-propagation scales linearly with the number of layers and the
number of samples, we define the total computational cost~$C$ as
follows:
\begin{align}
 C =\sum_{\ell=1}^r 2^{\ell-r} \,\sharp_\ell \, C_r,
 \label{eq:W_rmtr_mini_batch}
\end{align}
where the scaling factor\footnote{Uniform coarsening in 1D by a factor
of two is assumed.}~$2^{\ell-r}$ accounts for the difference between
the cost associated with level~$\ell$ and level~$r$.  The
symbol~$\sharp_\ell$ describes a number of gradient evaluations
performed on a level~$\ell$ using full dataset $\pazocal{D}$.  For
instance, if we evaluate gradient three times on level~$\ell$
using~$n_b$ samples, then~$\sharp_\ell=3{n_b}/{|\pazocal{D}|}$.

All presented experiments were obtained using XC50 compute nodes
(Intel Xeon E5-2690 v3 processor, NVIDIA Tesla P100 graphics card) of
the Piz Daint supercomputer from the Swiss National Supercomputing
Centre~(CSCS).

\subsection{Numerical examples}

We investigate the convergence properties and the efficacy of the
proposed \malg\ algorithms using three numerical examples from the
field of classification and regression. 

\subsubsection{Hyperspectral image segmentation using Indian Pines dataset} 
\label{sec:IndianPines}

Our first example arises from soil segmentation using hyperspectral
images provided by the Indian Pines dataset~\cite{baumgardner2015220}.
The input data was gathered by an Airborne Visible/Infrared Imaging
Spectrometer (AVIRIS) sensor and consist of $ 145\times 145$ pixels
with $200$ spectral bands in the range from $400$ to $2500$ nm.  Out
of all pixels, only~$10,249$ contain labeled data, which we split
into~$7,174$ for training and $3,075$ for validation.  Each pixel is
assigned to one of the $16$ classes, representing the type of land
cover, e.g., corn, soybean, etc.  Segmentation is performed using
ResNet with dense residual blocks of width~$50$, and the \textit{ReLu}
activation function $\sigma$.  Parameters $K, T, \beta_1, \beta_2$ are
set to $K=T=3$ and $\beta_1=\beta_2= 10^{-3}$.  Moreover, we employ
the softmax hypothesis function together with the cross-entropy loss
function, defined as~${\mathcal{h}(\hat{c}_s, c_s) := c_s^T \log
  (\hat{c}_s)}$.  To train ResNet, we run variants of \malg,
\ASTR\ and \al{SGD} methods without momentum.  All methods are
terminated as soon as $\text{acc}_{\text{train}} > 0.98$ or
$\text{acc}_{\text{val}} > 0.98$.  Termination also occurs as soon as
$\sum_{i=1}^{15}(\text{acc}_{\text{train}})_e-(\text{acc}_{\text{train}})_{e-i}<0.001$
or $\sum_{i=1}^{15}(\text{acc}_{\text{val}})_e - (\text{acc}_{\text{val}})_{e-i} <0.001$,
where $(\text{acc}_{\text{train/val}})_e$ is defined as the
ratio between the number of correctly classified samples from the
train/validation dataset and the total number of samples in the
train/validation dataset for a given epoch~$e$.

\begin{figure}[ht]
  \centering
\includegraphics{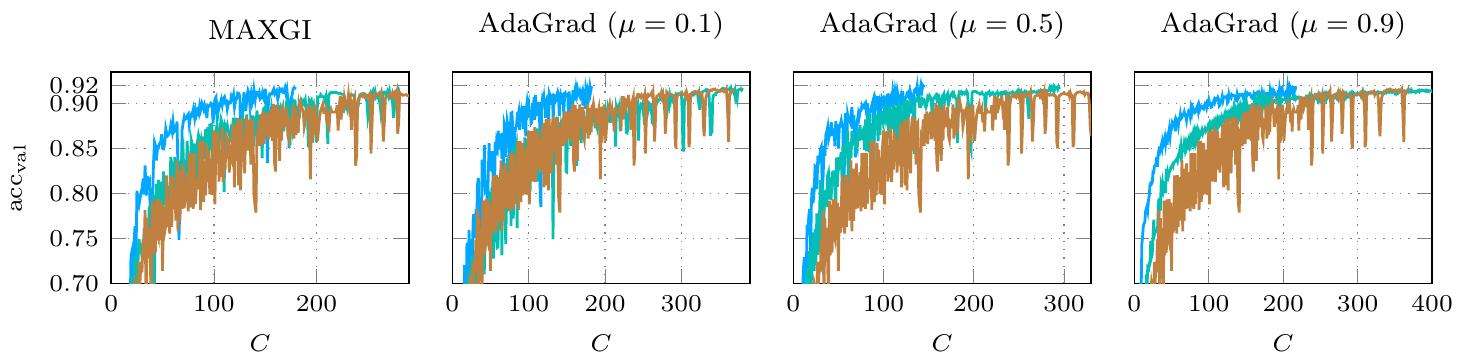}  
\includegraphics{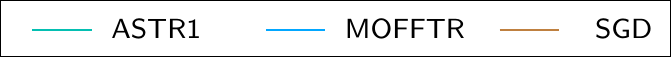}  
\caption{Indian Pines example: The validation accuracy as a
function of the total computational cost $C$.
The run with the highest validation accuracy is reported among $10$ independent runs. 
The experiments were performed using a batch size of $1,024$.} 
  \label{fig:indian_pines}
\end{figure}

\begin{table}[ht]
\centering
\begin{tabular}{|l|l||ll||ll||l|l||l|l|}
\hline
\multicolumn{2}{|c||}{\multirow{3}{*}{{Method} } }& 	\multicolumn{8}{c|}{{Batch-size (batch-size/$|\pazocal{D}|$)}}\\  \cline{3-10}
\multicolumn{2}{|c||}{ }	& \multicolumn{2}{c||}{{512 (6\%)}}                                  & \multicolumn{2}{c||}{{1,024 (12\%)}}                                        & \multicolumn{2}{c||}{{4,096 (50\%)}}    & \multicolumn{2}{c|}{{8,100 (100\%)}}                                 \\ \cline{3-10} 
\multicolumn{2}{|c||}{ } & \multicolumn{1}{c|}{${C}$} & $\text{{acc}}_\text{{val}}$ & \multicolumn{1}{c|}{${C}$} & $\text{{acc}}_\text{{val}}$ & \multicolumn{1}{c|}{${C}$} & $\text{{acc}}_\text{{val}}$ & \multicolumn{1}{c|}{${C}$} & $\text{{acc}}_\text{{val}}$  \\ \hline \hline
 \multicolumn{1}{|c}{ }&\al{SGD}         & \multicolumn{1}{r|}{359}                 & 92.5\%                  & \multicolumn{1}{r|}{369}                 & 91.7\%                  & \multicolumn{1}{r|}{1,230}                & 92.2\%        & {3,093} &   $92.1\%$      \\ \hline \hline

\multirow{4}{*}{\rotatebox[origin=c]{90}{{\ASTR}}} & \al{AdaGrad($\mu=0.1$)}           & \multicolumn{1}{r|}{230}             & 92.6\%                  & \multicolumn{1}{r|}{379}                 & 91.8\%                  & \multicolumn{1}{r|}{979}               & 92.3\% & 2,748 & 92.4\%               \\ \cline{2-10} 
                                 & \al{AdaGrad($\mu=0.5$)}         & \multicolumn{1}{r|}{143}                 & 92.4\%                  			& \multicolumn{1}{r|}{295}                 & 91.8\%                  & \multicolumn{1}{r|}{853}                 & 92.3\%     &     3,171 & 91.9\%         \\ \cline{2-10} 
                                 & \al{AdaGrad($\mu=0.9$)}         & \multicolumn{1}{r|}{232}                 & 91.4\%                  			& \multicolumn{1}{r|}{404}                 & 91.6\%                  & \multicolumn{1}{r|}{976}                & 92.0\%      &     2,626 & 91.7\%          \\ \cline{2-10} 
                                 & \al{MAXGI}                      & \multicolumn{1}{r|}{195}                         & 92.4\%                  			& \multicolumn{1}{r|}{281}                 & 91.7\%                  & \multicolumn{1}{r|}{884}                 & 92.2\%      &     2,719 & 91.8\%           \\ \hline \hline

\multirow{4}{*}{\rotatebox[origin=c]{90}{{\malg}}} 
				& \al{AdaGrad($\mu=0.1$)}      & \multicolumn{1}{r|}{125}                 & 92.4\%                  & \multicolumn{1}{r|}{183}                 & 91.7\%                  & \multicolumn{1}{r|}{422}                 & 92.3\%    & 1,683 &  $92.3\%$              \\ \cline{2-10} 
                                 & \al{AdaGrad($\mu=0.5$)}      	& \multicolumn{1}{r|}{113}                 & 92.1\%                  & \multicolumn{1}{r|}{145}              & 92.1\%                  & \multicolumn{1}{r|}{412}                 & 92.2\%    & 1,394 &  $92.3\%$                 \\ \cline{2-10} 
                                 & \al{AdaGrad($\mu=0.9$)}      	& \multicolumn{1}{r|}{118}                & 91.5\%                  & \multicolumn{1}{r|}{217}               & 91.9\%                  & \multicolumn{1}{r|}{453}                 & 92.1\%    & 1,760 &  $91.9\%$                \\ \cline{2-10} 
                                 & \al{MAXGI}                   		& \multicolumn{1}{r|}{96}                  & 92.3\%                   & \multicolumn{1}{r|}{174}               & 91.7\%                  & \multicolumn{1}{r|}{334}                 & 91.9\%   & 1,292 &  $92.4\%$                 \\ \hline
\end{tabular}
\caption{Indian Pines example: The total computational cost~$C$ and
  validation accuracy~$\text{acc}_\text{val}$ required to train
  ResNet. The best result in terms of validation accuracy is reported
  among $10$ independent runs.} 
\label{tab:indian_pines}
\end{table}

Table~\ref{tab:indian_pines} reports the total computational cost and
validation accuracy required by all solution strategies, for
increasing noise (i.e.,~decreasing batch size).  As can be seen in
this table, all solution strategies achieve comparable validation
accuracy for a given batch size.  However, the computational cost of
all variants of \malg\ is smaller than that of their single-level
counterparts and of the \al{SGD} method, see also
Figure~\ref{fig:indian_pines}.  Compared to \ASTR, the speedup factor
fluctuates from $1.2$ to $2.6$.  Compared to \al{SGD} method, the
speedup is higher as it ranges from~$1.7$ to $3.7$.  Moreover, our
results suggest that the speedup can be consistently observed even for
small batch sizes.  This empirically confirms that the proposed
algorithmic framework retains enhanced convergence of the multilevel
methods and at the same time is insensitive to the subsampling noise
as a majority of OFFO methods.

\subsubsection{Surrogate modelling of parametric neutron diffusion-reaction (NDR)}
\label{sec:ndr}

\begin{figure}[ht]
 \centering 
 \includegraphics[scale=0.205]{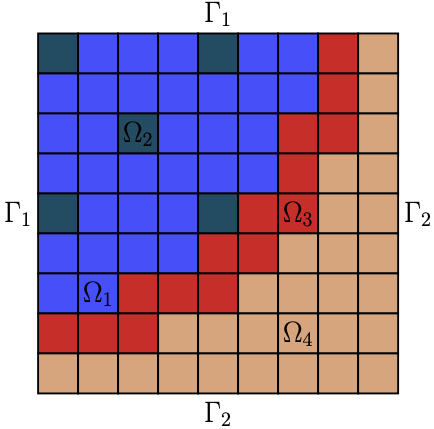}
 \quad
 \includegraphics[scale=0.0725]{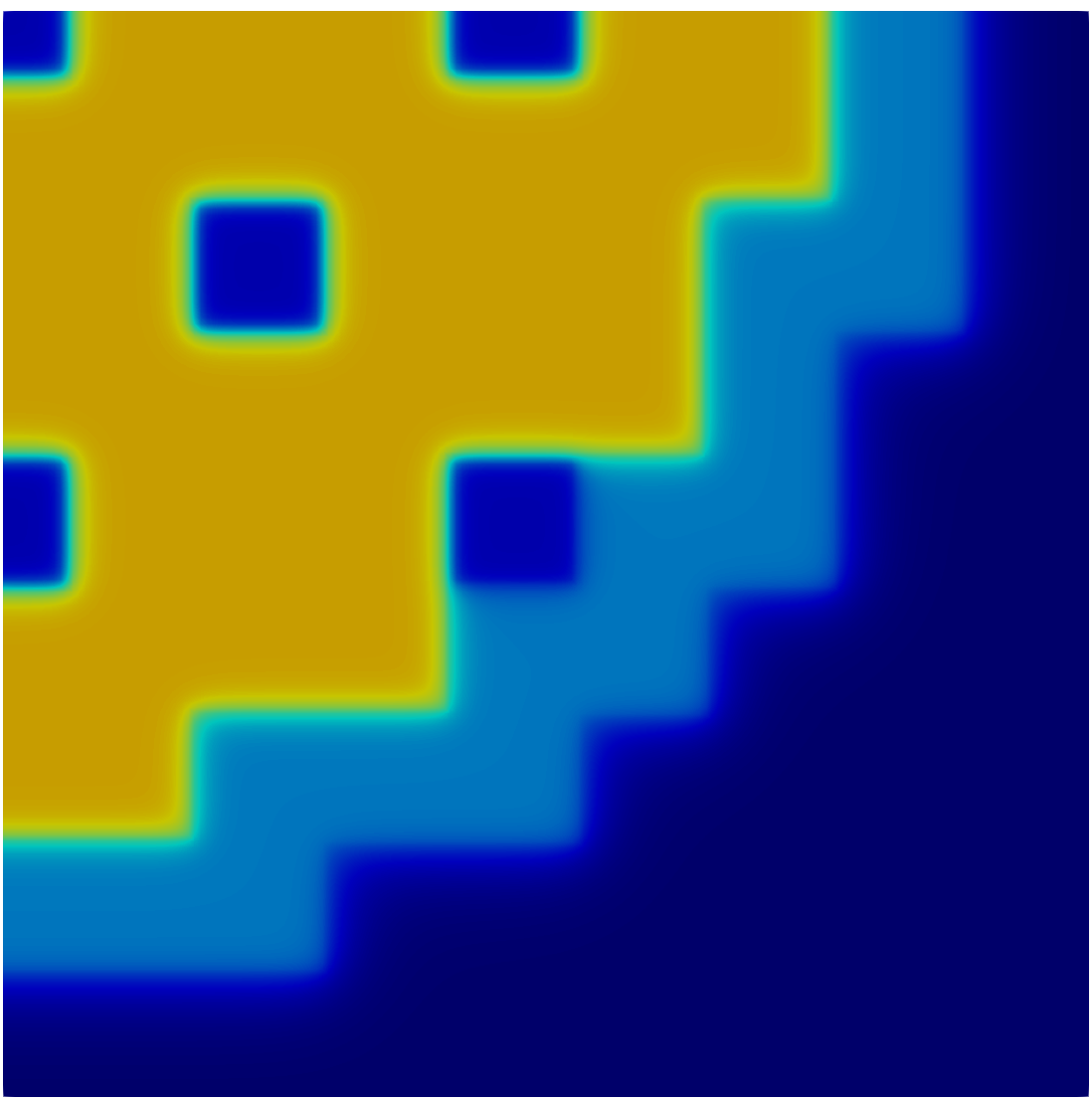}
 \quad
 \includegraphics[scale=0.0725]{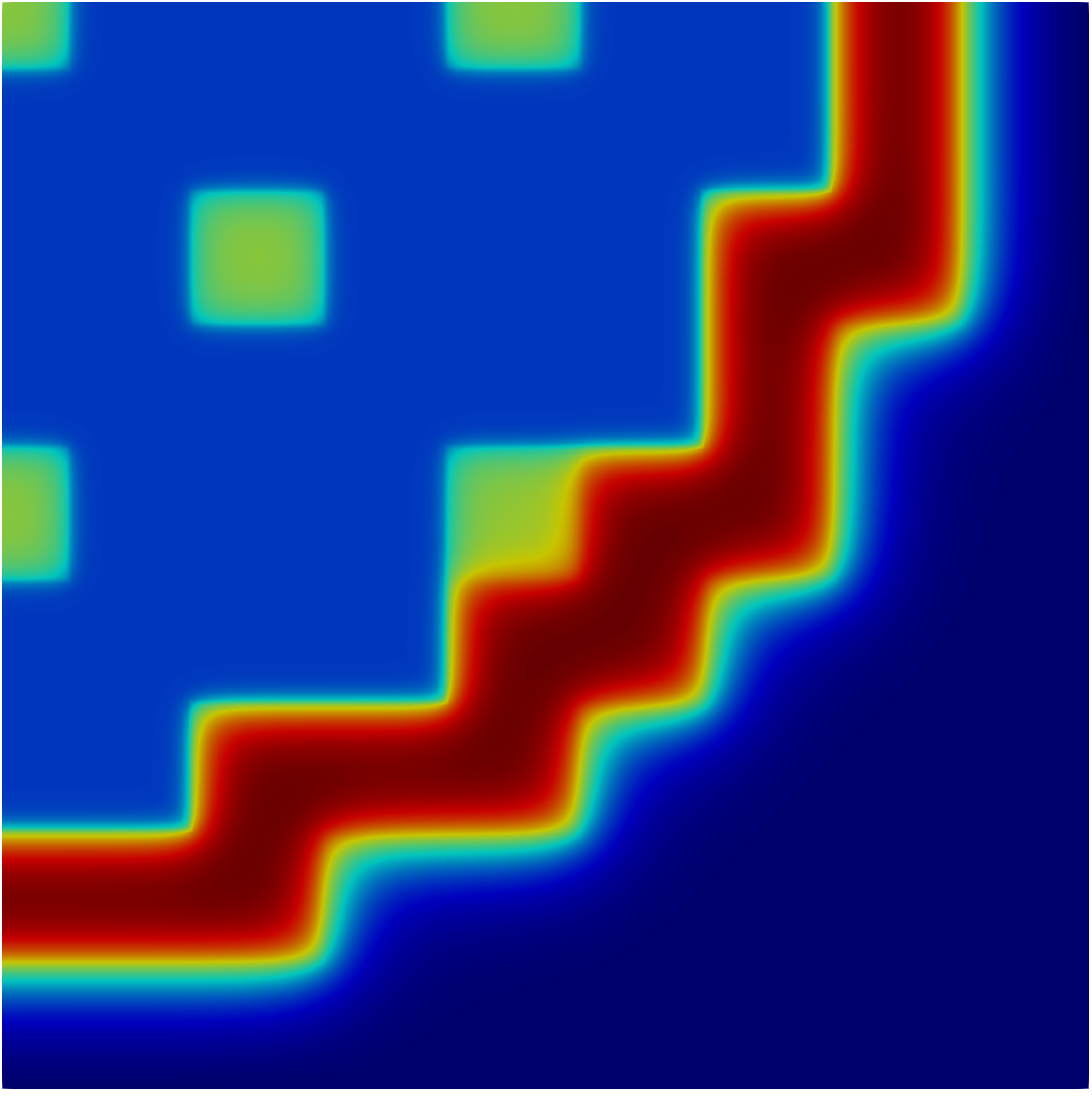}
 \quad
 \includegraphics[scale=0.0725]{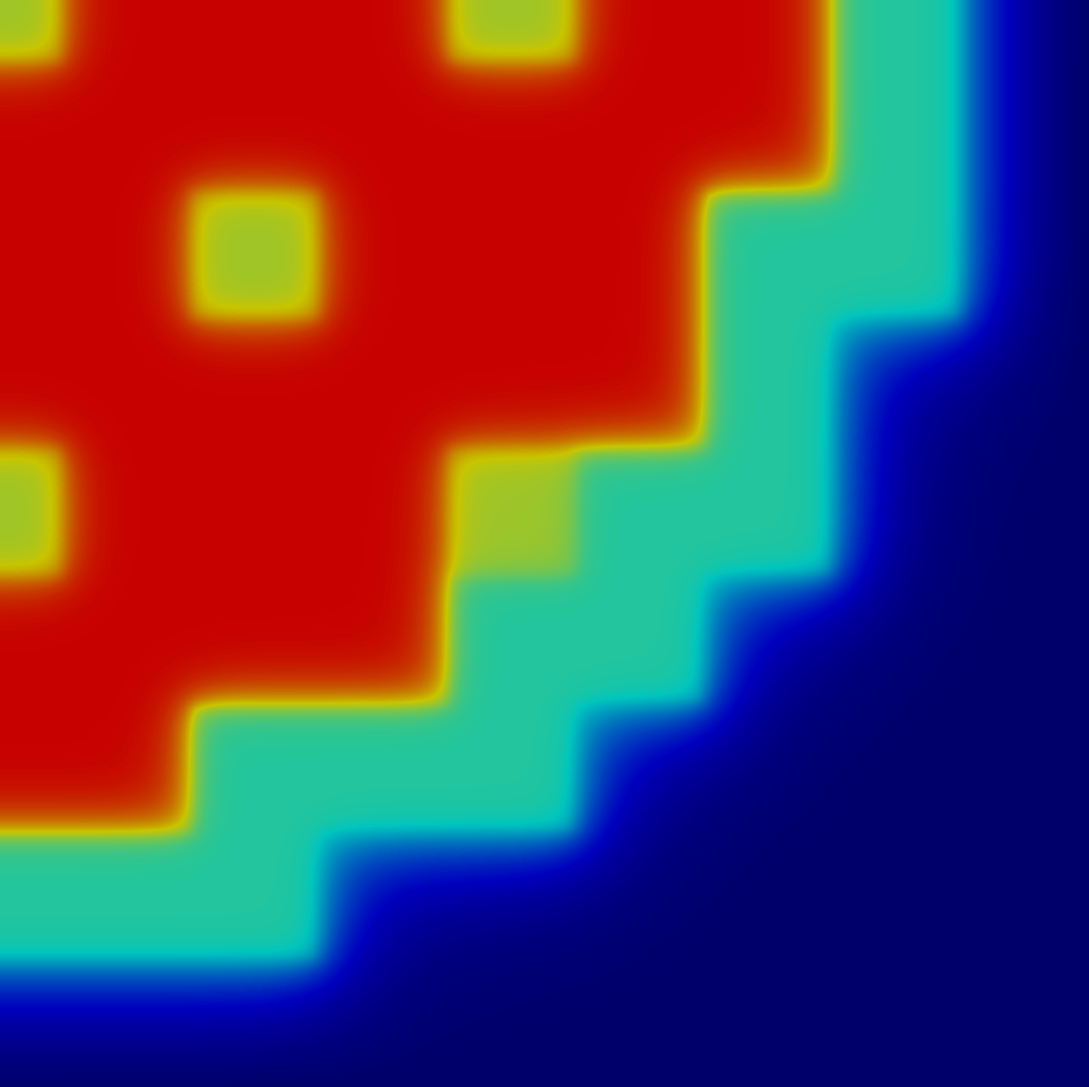}
 \quad
 \includegraphics[scale=0.145]{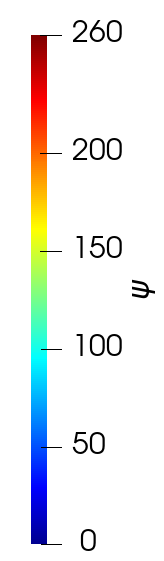}
 \caption{Left: Computational domain used for the creation of the NDR
   dataset. Each subdomain is illustrated by a different color. 
 Middle-Right: Example of samples contained in the NDR dataset.}
 \label{fig:datasets_NDR}
 \end{figure}

Our second example considers the construction of a surrogate model for
a parametric neutron diffusion-reaction problem with spatially-varying
coefficients and an external source.  The goal is to construct a
surrogate that can predict the average neutron flux for a given set of
parameters.  To this aim, we generate a dataset of $3,000$ samples,
which we split into~$2,600$ samples for training and $400$ for
testing.  Following~\cite{PrinRagu19}, the computational
domain~$\Omega=(0, 170)^2$ is heterogeneous and consists of four
different material regions, denoted by $\Omega_1, \ldots, \Omega_4 $
(see the left panel of Figure~\ref{fig:datasets_NDR}).

 The strong form of the problem is given as 
\begin{align}
 \nabla \cdot [ D(x) \nabla \psi(x) ] + \alpha (x) \nabla \psi(x) &= q(x), \qquad &&\text{in} \ \Omega, \nonumber \\ 
\psi(x) &= 0, &&\text{on} \ \Gamma_{1} := [0, 170] \times \{1\} \cup \{0\} \times [0, 170] , \label{eq:neutron}\\ 
 D(x) \nabla \psi(x) \cdot n(x) &= 0, &&\text{on} \ \Gamma_{2} := [0, 170] \times \{0\} \cup \{1\} \times [0, 170] , \nonumber
\end{align}
where~$x$ denotes the spatial coordinates and $\psi$ is the neutron flux from $\Omega$ to $\Re$.
Functions~$D, \alpha, q$ are defined as
${D(x) = \sum_{i=1}^4 \mathbbm{1}_{\Omega_i}(x) D_i}$,
${q(x) = \sum_{i=1}^3 \mathbbm{1}_{\Omega_i}(x) q_i}$,
and ${\alpha(x) = \sum_{i=1}^4 \mathbbm{1}_{\Omega_i}(x) \alpha_i}$, respectively
($\mathbbm{1}_{\Omega_i}$ denotes the indicator function of the domain
$\Omega$).  Problem~\eqref{eq:neutron} is then parametrized using $11$ parameters,
which we sample from a uniform distribution~$\pazocal{U}(a,b)$,
specified by lower (a) and upper (b) bounds.  
More precisely, diffusion coefficients $\{D_i\}_{i=1}^3$ are sampled
from $\pazocal{U}(0.15, 0.6)$, while~$D_4$ is sampled from~$\pazocal{U}(0.2, 0.8)$. 
Reaction coefficients~$\alpha_{1}, \ldots, \alpha_{4}$ take on values
from $\pazocal{U}(0.0425, 0.17)$, $\pazocal{U}(0.065, 0.26)$, $\pazocal{U}(0.04, 0.16)$, $\pazocal{U}(0.005, 0.02)$, respectively. 
The values of sources~$\{ q_i \}_{i=1}^3$ are sampled
from~$\pazocal{U}(5, 20)$, while the value value of $q_4$
is set to $0$.  For each set of parameters/input features, we create a
target $c_s = \int_{\Omega} \psi(x) \ dx /\int_{\Omega} dx$ by
solving~\eqref{eq:neutron} numerically using the finite element method
with a quadrilateral mesh ($500$ elements in both spatial dimensions).  

To build the desired surrogate, we train ResNet with dense residual
blocks of width $10$, \textit{tanh} activation function $\sigma$ and
parameters~$T=K=3$ and~$\beta_1=\beta_2=10^{-4}$.  
The identity hypothesis function and mean square loss functional are
defined as~${\mathcal{h}(\hat{c}_s, c_s) := \| c_s - \hat{c}_s \|_2^2
}$.

As common for regression problems, we train ResNet using variants of
the \malg, \ASTR, \al{SGD} and \al{ADAM} methods with momentum.  
The details of how to handle the momentum in the multilevel framework
can be found in~\cite[Appendix A]{kopanicakova2022globally}.  
The training is performed for a fixed computational budget, i.e., all 
solution strategies terminate as soon as~${C > C_{\text{max}}}$, where
we set~$C_{\text{max}}$ to $2,000$. 

\begin{table}[ht]
\centering
\begin{tabular}{|ll||r|r||r|r|r|}
\hline
\multicolumn{2}{|c||}{\multirow{3}{*}{{Method} } }& 	\multicolumn{4}{c|}{{Batch-size (batch-size/$|\pazocal{D}|$)}}\\  \cline{3-6}
                                  &         				& \multicolumn{2}{c||}{{256 (25\%)}}                              & \multicolumn{2}{c|}{{1,024 (100\%)}}                                                                \\ \cline{3-6} 
\multicolumn{2}{|c||}{}                                         & ${f_{\text{{train}}} (\times 10^3)}$ & 	${f_{\text{{val}}} (\times 10^3)}$ & ${f_{\text{{train}}} (\times 10^3)}$ & ${f_{\text{{val}}} (\times 10^3)}$ \\ \hline \hline
 \multicolumn{1}{|c}{ }                                        & \al{SGD}   & $0.61 \pm 0.30$              & $3.57 \pm 1.44$                    & $1.87 \pm 1.32$             & $6.87 \pm 1.09$                                 \\ \cline{1-6} 
 \multicolumn{1}{|c}{ }                                        & \al{ADAM}   &  $0.39 \pm 0.23$             & $5.77 \pm 0.79$                   & $0.47 \pm 0.50$               & $5.35 \pm 0.23$                                 \\ \hline \hline

\multirow{4}{*}{\rotatebox[origin=c]{90}{{\ASTR}}}            &\multicolumn{1}{|l||}{   \al{AdaGrad($\mu=0.1$)} }& $0.56 \pm 0.13$      &       $5.41 \pm 0.72$                                          &       $1.01 \pm 0.37$                    	&           $5.72 \pm 1.18$                                     \\ \cline{2-6} 
\multicolumn{1}{|l|}{}                                                 & {\al{AdaGrad($\mu=0.5$)}} 	&    $0.47 \pm 0.32$           	&         $4.62 \pm 2.71$                                       &         $0.79 \pm 0.12$                      &    $7.04 \pm 4.93$                                            \\ \cline{2-6} 
\multicolumn{1}{|l|}{}                                                 & {\al{AdaGrad($\mu=0.9$)}} 	&    $0.54 \pm 0.17$           	&           $7.43 \pm 2.48$                                     &         $0.86 \pm 0.10$                       &     $4.93 \pm 2.56$                                           \\ \cline{2-6} 
\multicolumn{1}{|l|}{}                                                 & \al{MAXGI}              		&      $0.53 \pm 0.35$       		&           $5.67 \pm 1.23$                                     &          $1.04 \pm 0.29$                  	&      $4.98 \pm 1.28$                                      \\  \hline \hline

\multirow{4}{*}{\rotatebox[origin=c]{90}{{\malg}}}&\multicolumn{1}{|l||}{  \al{AdaGrad($\mu=0.1$)}} & $0.48 \pm 0.29$            & $0.89 \pm 0.31$                     & $0.48 \pm 0.73$             & $1.02 \pm 0.42$                                  \\ \cline{2-6} 
\multicolumn{1}{|l|}{}                                                 & \al{AdaGrad($\mu=0.5$)}	 & $0.39 \pm 0.31$            & $0.95 \pm 0.29$                    & $0.63 \pm 0.45$              & $1.29 \pm 0.53$                                  \\ \cline{2-6} 
\multicolumn{1}{|l|}{}                                                 & \al{AdaGrad($\mu=0.9$)} 	& $0.45 \pm 0.25$             & $1.06 \pm 0.36$                      & $0.62 \pm 0.41$             & $1.36 \pm 0.49$                                 \\ \cline{2-6} 
\multicolumn{1}{|l|}{}                                                 & \al{MAXGI}              		&   $0.51 \pm 0.14$           & $0.91 \pm 0.30$                                    		&$0.73 \pm 0.47$             & $1.38 \pm 0.45$                                 \\  \hline
\end{tabular}
\caption{NDR example: The mean and standard deviation of the value of train ($f_{\text{train}}$) and validation ($f_{\text{val}}$) loss function achieved by training ResNet for $2,000$ epochs.
The statistics is obtained from $10$ independent runs.}
\label{table:NDR}
\end{table}

\begin{figure}[ht]
\centering
\includegraphics{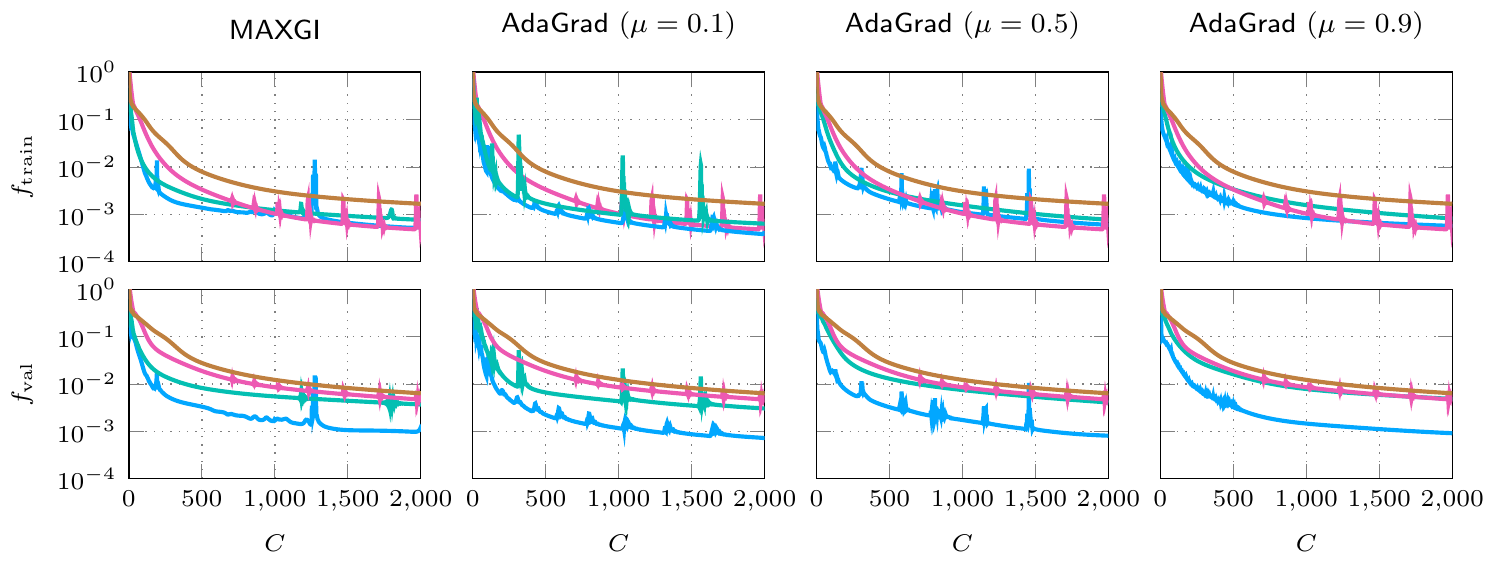}
\includegraphics{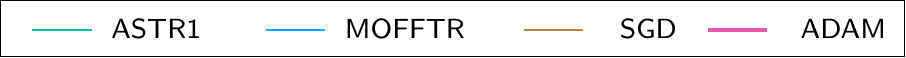}    
  \caption{NDR example: The history of train ($f_{\text{train}}$) and
    validation ($f_{\text{val}}$) loss function obtained while
    training ResNet for $2,000$ epochs with batch size $1,024$ (full
    dataset). The result with lowest $f_{\text{val}}$ among $10$ independent runs is reported.}
  \label{fig:neutron}
\end{figure}

We observe in Table~\ref{table:NDR} that all solution strategies,
except \al{SGD}, achieve comparable values of the training loss
($f_{\text{train}}$). Interestingly, we also notice that all variants of the 
\malg\ algorithm generalize better, i.e., reach a lower value of the
validation loss ($f_{\text{val}}$) than all single-level methods.  
This phenomenon can be also observed in Figure~\ref{fig:neutron}. 
The results presented in Table~\ref{table:NDR} suggest that the
\malg\ methods preserve good generalization properties 
in the presence of noise, a property of particular interest in
surrogate modeling and other scientific applications that require
reliable solutions. 

\subsubsection{Image classification using the SVNH dataset}

Our last numerical example is associated with an image classification
using SVNH dataset~\cite{netzer2011reading}.  Each image is
represented by~${32 \times 32}$ pixels and contains overlapping digits
from $0$ to $9$.  This dataset consists of $99,289$ samples, from
which $73,257$ are used for training and $26,032$ for testing
purposes.  We pre-process all samples by standardizing the images, so
that pixel values lie in the range~$[0, 1]$, and by subtracting the
mean from each pixel.  In addition, we use standard data augmentation
techniques, in particular image rotation, horizontal and vertical
shift, and horizontal flip.  The image classification is performed
using ResNet with convolutional residual blocks (32 filters),
\textit{ReLu} activation function and parameters~$T = 3$,~$K=5$,
${\beta_1=10^{-3}}$, ${\beta_2=0.005}$.  Moreover, we use the softmax
hypothesis function with the cross-entropy loss function.

The training of ResNets is performed using a batch size of~$512$ and
the same stopping criterion as that used for the Indian Pines example
in Section~\ref{sec:IndianPines}.  For this experiment, we consider
the \malg, \ASTR\ and \al{SGD} methods without momentum.
Table~\ref{table:svnh} reports the computational cost and validation
accuracy achieved by all methods.  The study is performed with respect
to an increasing number of refinement levels\footnote{The number of
levels utilized by the \malg\ algorithms increases linearly with
refinement level.}.

\begin{figure}[ht]
  \centering
\includegraphics{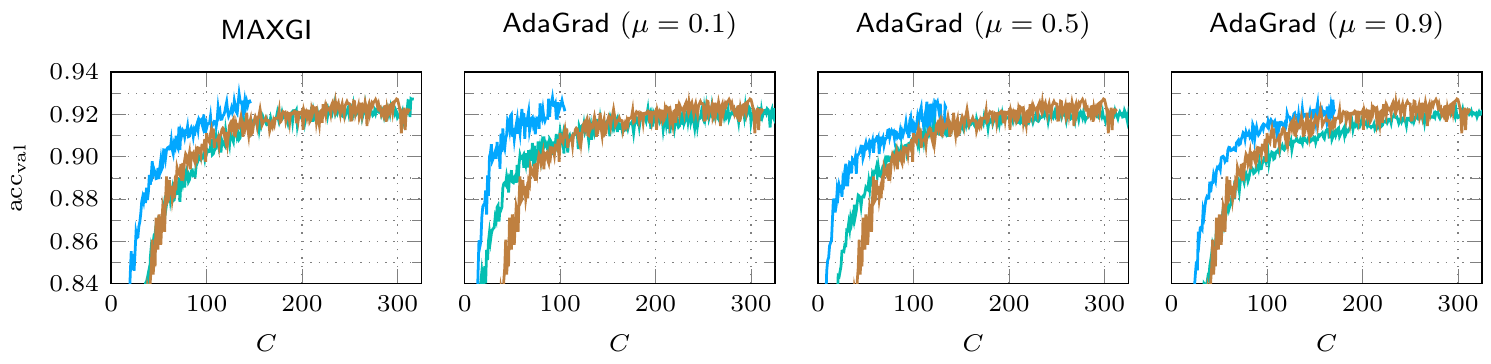}  
\includegraphics{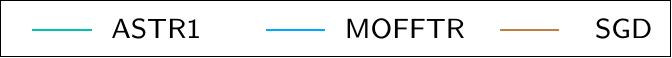}
  \caption{SVNH example: The validation accuracy as a function of
    total computational cost obtained while training ResNet with $17$
    residual block (three refinement levels). The run with the highest
    validation accuracy is selected among $10$ independent runs.}  
  \label{fig:svnh}
\end{figure}

\begin{table}[ht]
\centering
\begin{tabular}{|l|l||l|l||l|l||l|l|}
\hline
\multicolumn{2}{|c||}{\multirow{3}{*}{{Method} } }& 	\multicolumn{6}{c|}{{Levels (residual blocks)}}\\  \cline{3-8}
\multicolumn{2}{|c||}{ }	& \multicolumn{2}{c||}{{2 (9)}}                                  & \multicolumn{2}{c||}{{3 (17)}}                                        & \multicolumn{2}{c|}{{4 (33)}}         \\ \cline{3-8} 
\multicolumn{2}{|c||}{ } & \multicolumn{1}{c|}{${C}$} & $\text{{acc}}_\text{{val}}$ & \multicolumn{1}{c|}{${C}$} & $\text{{acc}}_\text{{val}}$ & \multicolumn{1}{c|}{${C}$} & $\text{{acc}}_\text{{val}}$  \\ \hline \hline
 \multicolumn{1}{|c}{ }&\al{SGD}           & \multicolumn{1}{r|}{280}                 & 92.4\%                  & \multicolumn{1}{r|}{313}                & 92.7\%     &   336    & 92.9\%          \\ \hline \hline

\multirow{4}{*}{\rotatebox[origin=c]{90}{{\ASTR}}} & {\al{AdaGrad($\mu=0.1$)}} 		&     325     &      92.4\%        	&      $361$ 		&    92.6\%        	&     341       	&   92.8\%     \\ \cline{2-8} 
                                 & \al{AdaGrad($\mu=0.5$)}         	&     334     &      92.3\%       	&       $333$		&    92.4\%          	&     328     	&   92.7\%         \\ \cline{2-8} 
                                 & \al{AdaGrad($\mu=0.9$)}         	&      351    &    92.1\%          	&       $358$		&    92.4\%         	&     369     	&  92.6\%            \\ \cline{2-8} 
                                 & \al{MAXGI}                      		&  313        &    92.3\%          	&   $317$    		&   92.7\%         	&     325      	&   92.8\%          \\ \hline \hline

\multirow{4}{*}{\rotatebox[origin=c]{90}{{\malg}}} 
				& \al{AdaGrad($\mu=0.1$)}     	&     134   	  &     92.3\%         	&   $105$    		&       92.7\%       	&     139      & 		92.9\% 		\\ \cline{2-8} 
                                 & \al{AdaGrad($\mu=0.5$)}      	&     148  	  &     92.2\%         	&   $135$  		&       92.6\%      	&      144     &          92.8\%       	\\ \cline{2-8} 
                                 & \al{AdaGrad($\mu=0.9$)}      	&     169  	  &     92.1\%        	&   $171$   		&      92.6\%       	&       151     &          92.8\%        		\\ \cline{2-8} 
                                 & \al{MAXGI}                   		&     158     &     92.5\%          	&  $146$       		&      92.9\%          	&       143	    &         93.1\%     		\\ \hline
\end{tabular}
\caption{SVNH example: The total computational cost~$C$ and validation
  accuracy~$\text{acc}_\text{val}$, which were required to train
  ResNet with increasing depth (refinement level). 
  The best result in terms of validation accuracy is reported among $10$
  independent runs.} 
\label{table:svnh}
\end{table}

These results suggest that increasing the network depth enhances its
representation capacity, which in turn justifies a higher
computational cost.  We also notice that all variants of the
\malg\ method require lower computational cost than single-level
solution strategies while achieving comparable accuracy.  The obtained
computational gains vary from the factor of $1.7$ to $3.4$.  This can
be observed independently of the refinement level, which indicates
that the \malg\ algorithm has a clear potential to speed up the
training of large-scale networks, such as the ones used in real-life
applications.  We also note that, among all variants of the
\malg\ framework, the choice of the \al{MAXGI} weights yields higher
validation accuracy than that obtained with the AdaGrad-like variants.
We finally observe that, among all AdaGrad-like variants, the
configuration with~$\mu=0.1$ achieves the highest validation accuracy
as well as the lowest computational cost.  The fact that a similar
observation can be made also for the corresponding single-level
\ASTR\ methods suggests that exploring values of~$\mu$ other than the
traditional $\mu=0.5$ might be beneficial in practice.

\section{Conclusion}\label{section:conclusion}

We have presented a class of multilevel algorithms for unconstrained
minimization which do not require the computation of the objective
function's value. We have also shown that the performance of
algorithms of the class is competitive in the presence of
noise-induced by subsampling in the context of deep neural network
training. The convergence theory of two interesting subclasses has
been analyzed and shown to match the state of the art.  Our
experiments also indicate that, although currently not covered by our
theory, the benefits of the multilevel approach are preserved when
momentum is added to the framework. The authors are well aware that
only continued experimentation will reveal the true practical value of
the present proposal, but they note that the first numerical tests are
encouraging.

{\footnotesize

\section*{\footnotesize Acknowledgements}
A.~Kopani\v{c}\'akov\'a gratefully acknowledges support of the Swiss
National Science Foundation through the project "Multilevel training of
DeepONets — multiscale and multiphysics applications" (grant
no.~206745), and partial support of Platform for Advanced Scientific 
Computing (PASC) under the project EXATRAIN. 
Ph.~L.~Toint acknowledges the continued and friendly partial support of ANITI.


\begin{thebibliography}{10}

\bibitem{AlexLewi01}
N.~M. Alexandrov and R.~L. Lewis.
\newblock An overview of first-order model management for engineering
  optimization.
\newblock {\em Optimization and Engineering}, 2:413--430, 2001.

\bibitem{baumgardner2015220}
M.~F. Baumgardner, L.~L. Biehl, and D.~A. Landgrebe.
\newblock 220 band {A}{V}{I}{R}{I}{S} {H}yperspectral {I}mage {D}ata {S}et:
  {J}une 12, 1992 {I}ndian {P}ine {T}e{S}ite 3.
\newblock Purdue University Research Repository, 10(7):991, 2015.

\bibitem{BeckHall18}
A.~Beck and N.~Hallak.
\newblock Optimization problems involving group sparsity terms.
\newblock {\em Mathematical Programming}, 2018.
\newblock online.

\bibitem{BorzKuni06}
A.~Borzi and K.~Kunisch.
\newblock A globalisation strategy for the multigrid solution of elliptic
  optimal control problems.
\newblock {\em Optimization Methods and Software}, 21(3):445--459, 2006.

\bibitem{BouaSchn98}
A.~Bouaricha and R.~B. Schnabel.
\newblock Tensor methods for large sparse systems of nonlinear equations.
\newblock {\em Mathematical Programming}, 82(3):377--412, 1998.

\bibitem{BrigHensMcCo00}
W.~L. Briggs, V.~E. Henson, and S.~F. McCormick.
\newblock {\em A Multigrid Tutorial}.
\newblock SIAM, Philadelphia, USA, 2nd edition, 2000.

\bibitem{CartGoulToin22}
C.~Cartis, N.~I.~M. Gould, and Ph.~L. Toint.
\newblock {\em Evaluation complexity of algorithms for nonconvex optimization}.
\newblock Number~30 in MOS-SIAM Series on Optimization. SIAM, Philadelphia,
  USA, June 2022.

\bibitem{chang2017multi}
B.~Chang, L.~Meng, E.~Haber, F.~Tung, and D.~Begert.
\newblock Multi-level residual networks from dynamical systems view.
\newblock arXiv:1710.10348, 2017.

\bibitem{Chat13}
I.~Chatzigeorgiou.
\newblock Bounds on the {Lambert} function and their application to the outage
  analysis of user cooperation.
\newblock {\em {IEEE} Communications Letters}, 17(8):1505--1508, 2013.

\bibitem{ConnGoulToin92}
A.~R. Conn, N.~I.~M. Gould, and Ph.~L. Toint.
\newblock {\em {\sf LANCELOT}: a {F}ortran package for large-scale nonlinear
  optimization ({R}elease {A})}.
\newblock Number~17 in Springer Series in Computational Mathematics. Springer
  Verlag, Heidelberg, Berlin, New York, 1992.

\bibitem{ConnGoulToin00}
A.~R. Conn, N.~I.~M. Gould, and Ph.~L. Toint.
\newblock {\em Trust-Region Methods}.
\newblock Number~1 in MOS-SIAM Optimization Series. SIAM, Philadelphia, USA,
  2000.

\bibitem{Corletal96}
R.~M. Corless, G.~H. Gonnet, D.~E. Hare, D.~J. Jeffrey, and D.~E. Knuth.
\newblock On the {L}ambert {W} function.
\newblock {\em Advances in Computational Mathematics}, 5:329--–359, 1996.

\bibitem{cyr2019multilevel}
E.~C. Cyr, S.~G{\"u}nther, and J.~B. Schroder.
\newblock Multilevel initialization for layer-parallel deep neural network
  training.
\newblock {\em International Journal of Computing and Visualization in Science
  and Engineering}, 2021.

\bibitem{DuchHazaSing11}
J.~Duchi, E.~Hazan, and Y.~Singer.
\newblock Adaptive subgradient methods for online learning and stochastic
  optimization.
\newblock {\em Journal of Machine Learning Research}, 12, July 2011.

\bibitem{chollet2015keras}
F.~Chollet et~al.
\newblock Keras, 2015.
\newblock \url{https://keras.io}.

\bibitem{Fish98}
M.~Fisher.
\newblock Minimization algorithms for variational data assimilation.
\newblock In {\em Recent Developments in Numerical Methods for Atmospheric
  Modelling}, pages 364--385, Reading, UK, 1998. European Center for
  Medium-Range Weather Forecasts.

\bibitem{FujiKojiNaka97}
K.~Fujisawa, M.~Kojima, and K.~Nakata.
\newblock Exploiting sparsity in primal-dual interior-point methods for
  semidefinite programming.
\newblock {\em Mathematical Programming, Series~B}, 79(1--3):235--253, 1997.

\bibitem{Kopanicakova_2020c}
L.~Gaedke-Merzh{\"a}user, A.~{Kopani{\v{c}}{\'a}kov{\'a}}, and R.~Krause.
\newblock Multilevel minimization for deep residual networks.
\newblock In {\em ESAIM. Proceedings and Surveys.} 71:131-144, 2021.

\bibitem{Gay96}
D.~M. Gay.
\newblock Automatically finding and exploiting partially separable structure in
  nonlinear programming problems.
\newblock Technical report, Bell Laboratories, Murray Hill, New Jersey, USA,
  1996.

\bibitem{GelmMand90}
E.~Gelman and J.~Mandel.
\newblock On multilevel iterative methods for optimization problems.
\newblock {\em Mathematical Programming}, 48(1):1--17, 1990.

\bibitem{GoldWang93}
D.~Goldfarb and S.~Wang.
\newblock Partial-update {N}ewton methods for unary, factorable and partially
  separable optimization.
\newblock {\em SIAM Journal on Optimization}, 3(2):383--397, 1993.

\bibitem{GratJeraToin22c}
S.~Gratton, S.~Jerad, and Ph.~L. Toint.
\newblock Convergence properties of an objective-function-free optimization
  regularization algorithm, including an $\mathcal{O}(\epsilon^{-3/2})$
  complexity bound.
\newblock arXiv:2203.09947, 2022.

\bibitem{GratJeraToin22b}
S.~Gratton, S.~Jerad, and Ph.~L. Toint.
\newblock First-order objective-function-free optimization algorithms and their
  complexity.
\newblock arXiv:2203.01757, 2022.

\bibitem{GratJeraToin22a}
S.~Gratton, S.~Jerad, and Ph.~L. Toint.
\newblock Parametric complexity analysis for a class of first-order
  {A}dagrad-like algorithms.
\newblock arXiv:2203.01647, 2022.

\bibitem{GratMoufSartToinToma10}
S.~Gratton, M.~Mouffe, A.~Sartenaer, Ph.~L. Toint, and D.~Tomanos.
\newblock Numerical experience with a recursive trust-region method for
  multilevel nonlinear bound-constrained optimization.
\newblock {\em Optimization Methods and Software}, 25(3):359 -- 386, 2010.

\bibitem{GratSartToin08}
S.~Gratton, A.~Sartenaer, and Ph.~L. Toint.
\newblock Recursive trust-region methods for multiscale nonlinear optimization.
\newblock {\em SIAM Journal on Optimization}, 19(1):414--444, 2008.

\bibitem{GratToin09}
S.~Gratton and Ph.~L. Toint.
\newblock Approximate invariant subspaces and quasi-{N}ewton optimization
  methods.
\newblock {\em Optimization Methods and Software}, 25(4):507--529, 2010.

\bibitem{GrieToin82a}
A.~Griewank and Ph.~L. Toint.
\newblock On the unconstrained optimization of partially separable functions.
\newblock In M.~J.~D. Powell, editor, {\em Nonlinear Optimization 1981}, pages
  301--312, London, 1982. Academic Press.

\bibitem{haber2017stable}
E.~Haber and L.~Ruthotto.
\newblock Stable architectures for deep neural networks.
\newblock {\em Inverse Problems}, 34(1):014004, 2017.

\bibitem{haber2018learning}
E.~Haber, L.~Ruthotto, E.~Holtham, and S.-H. Jun.
\newblock Learning across scales---multiscale methods for convolution neural
  networks.
\newblock In {\em Thirty-Second AAAI Conference on Artificial Intelligence},
  2018.

\bibitem{he2016identity}
K.~He, X.~Zhang, S.~Ren, and J.~Sun.
\newblock Identity mappings in deep residual networks.
\newblock In {\em European conference on computer vision}, pages 630--645,
  Heidelberg, Berlin, New York, 2016. Springer Verlag.

\bibitem{ioffe2015batch}
S.~Ioffe and {Ch}. Szegedy.
\newblock Batch normalization: Accelerating deep network training by reducing
  internal covariate shift.
\newblock arXiv:1502.03167, 2015.

\bibitem{KingBa15}
D.~Kingma and J.~Ba.
\newblock {A}dam: A method for stochastic optimization.
\newblock In {\em Proceedings in the International Conference on Learning
  Representations (ICLR)}, 2015.

\bibitem{Kopanicakova2022e}
A.~Kopani\v{c}\'akov\'a.
\newblock On the use of hybrid coarse-level models in multilevel minimization methods.
\newblock arXiv:2211.15078, 2022.

\bibitem{kopanicakova2022globally}
A.~{Kopani{\v{c}}{\'a}kov{\'a}} and R.~Krause.
\newblock Globally convergent multilevel training of deep residual networks.
\newblock {\em SIAM Journal on Scientific Computing}, 0:S254--S280, 2022.

\bibitem{Korn97}
R.~Kornhuber.
\newblock Adaptive monotone multigrid methods for some non-smooth optimization
  problems.
\newblock In R.~Glowinski, J.~P\'{e}riaux, Z.~Shi, and O.~Widlund, editors,
  {\em Domain Decomposition Methods in Sciences and Engineering}, pages
  177--191. J. Wiley and Sons, Chichester, England, 1997.

\bibitem{Lasse05}
J.~B. Lasserre.
\newblock Convergent semidefinite relaxation in polynomial optimization with
  sparsity.
\newblock Technical report, LAAS-CNRS, 7, avenue du Colonel Roche, 31077
  Toulouse, France, November 2005.

\bibitem{LewiNash02}
M.~Lewis and S.~G. Nash.
\newblock Practical aspects of multiscale optimization methods for {VLSICAD}.
\newblock In Jason Cong and Joseph~R. Shinnerl, editors, {\em Multiscale
  Optimization and {VLSI/CAD}}, pages 265--291, Dordrecht, The Netherlands,
  2002. Kluwer Academic Publishers.

\bibitem{LewiNash04}
M.~Lewis and S.~G. Nash.
\newblock Model problems for the multigrid optimization of systems governed by
  differential equations.
\newblock {\em SIAM Journal on Scientific Computing}, 26(6):1811--1837, 2005.

\bibitem{MareRichTaka14}
J.~Mare\v{c}ek, P.~Richt\'arik, and M.~Tak\'a\v{c}.
\newblock Distributed block coordinate descent for minimizing partially
  separable functions.
\newblock Technical report, Department of Mathematics and Statistics,
  University of Edinburgh, Edinburgh, Scotland, 2014.

\bibitem{Nash00}
S.~G. Nash.
\newblock A multigrid approach to discretized optimization problems.
\newblock {\em Optimization Methods and Software}, 14:99--116, 2000.

\bibitem{netzer2011reading}
Y.~Netzer, T.~Wang, A.~Coates, A.~Bissacco, B.~Wu, and A.Y. Ng.
\newblock Reading digits in natural images with unsupervised feature learning,
  2011.

\bibitem{PrinRagu19}
Z.~M. Prince and J.~C. Ragusa.
\newblock Parametric uncertainty quantification using proper generalized
  decomposition applied to neutron diffusion.
\newblock {\em International Journal for Numerical Methods in Engineering},
  119(9):899–921, 2019.

\bibitem{queiruga2020continuous}
A.~F. Queiruga, N.~B. Erichson, D.~Taylor, and M.~W. Mahoney.
\newblock Continuous-in-depth neural networks.
\newblock arXiv:2008.02389, 2020.

\bibitem{ReddKaleKuma18}
S.~Reddi, S.~Kale, and S.~Kumar.
\newblock On the convergence of {A}dam and beyond.
\newblock In {\em Proceedings in the International Conference on Learning
  Representations (ICLR)}, 2018.

\bibitem{RobbMunr51}
H.~Robbins and S.~Monro.
\newblock A stochastic approximation method.
\newblock {\em The Annals of Mathematical Statistics}, 22(3):400--–407, 1951.

\bibitem{TielHint12}
T.~Tieleman and G.~Hinton.
\newblock Lecture~6.5-{RMSPROP}.
\newblock {COURSERA}: Neural Networks for Machine Learning, 2012.

\bibitem{Toin81a}
Ph.~L. Toint.
\newblock A note on sparsity exploiting quasi-{N}ewton methods.
\newblock {\em Mathematical Programming}, 21(2):172--181, 1981.

\bibitem{Toin81b}
Ph.~L. Toint.
\newblock Towards an efficient sparsity exploiting {N}ewton method for
  minimization.
\newblock In I.~S. Duff, editor, {\em Sparse Matrices and Their Uses}, pages
  57--88, London, 1981. Academic Press.

\bibitem{weinan2017proposal}
E.~Weinan.
\newblock A proposal on machine learning via dynamical systems.
\newblock {\em Communications in Mathematics and Statistics}, 5(1):1--11, 2017.

\bibitem{Yuan15}
Y.~Yuan.
\newblock Recent advances in trust region algorithms.
\newblock {\em Mathematical Programming, Series~A}, 151(1):249--281, 2015.

\bibitem{ZhouChenTangYangCaoGu20}
D.~Zhou, J.~Chen, Y.~Tang, Z.~Yang, Y.~Cao, and Q.~Gu.
\newblock On the convergence of adaptive gradient methods for nonconvex
  optimization.
\newblock arXiv:2080.05671, 2020.

\end{thebibliography}

}

\appendix

\section{Proof of Theorem~\ref{theorem:divs}}\label{proof-divs}

We first recall a useful technical result.

\llem{lemma:divsprop}{\cite[Lemma~3.5]{GratJeraToin22b} Consider and arbitrary $i\in \ii{n}$ and suppose
that there exists a $j_\varsigma$ such that
\beqn{termissmall}
\min\left[\frac{g_{i,j}^2}{\varsigma_i},\frac{g_{i,j}^2}{v_{i,j}}\right]
\leq \varsigma_i
\tim{ for } j \geq j_\varsigma.
\eeqn
Then
\beqn{termlbound}
\min\left[\frac{g_{i,j}^2}{\varsigma_i},\frac{g_{i,j}^2}{v_{i,j}}\right]
\geq \frac{g_{i,j}^2}{2\varsigma_i}
\tim{ for } j \geq j_\varsigma.
\eeqn
}

\paragraph*{Proof of Theorem~\ref{theorem:divs}}

We obtain from \req{fdecrease} and \req{fishr} that,
for $i \ge i_\vartheta$,
\[
f(x_{r,i_\vartheta})-f(x_{r,i+1})
\ge \frac{1}{2} \sum_{k=i_\vartheta+1}^i\sum_{j=1}^n\frac{g_{r,k,j}^2}{w_{r,k,j}}
    \left[\beta_{1,r}-\frac{\beta_{2,r}+\half\alpha^2 L}{w_{\ell,k,j}}\right].
\]
Moreover, \req{ws-divs} and \req{istar-divs} also ensure that, since
$i \ge i_\vartheta$,   
\[
\beta_{1,r}-\frac{\beta_{2,r}+\half\alpha^2L}{w_{\ell,k,j}} \ge \vartheta
\]
for all $j\in\ii{n}$, so that
\begin{align*}
f(x_{r,i_\vartheta}+1)-f(x_{r,i+1})
&\geq \frac{\vartheta}{2}\sum_{k=i_\vartheta+1}^i\sum_{j=1}^n\frac{g_{r,k,j}^2}{w_{r,k,j}}\\
&\geq \frac{\vartheta}{2} \sum_{k=i_\vartheta+1}^i\sum_{j=1}^n
\frac{g_{r,k,j}^2}{\max[\varsigma_j,v_{r,k}](k+1)^\mu}\\
&\geq \frac{\vartheta}{2(i+1)^\mu} \sum_{k=i_\vartheta+1}^i\sum_{j=1}^n
\min\left[\frac{g_{r,k,j}^2}{\varsigma_j},\frac{g_{r,k,j}^2}{v_{r,k}}\right]\\
&\geq \frac{\vartheta(i-i_\vartheta)}{2(i+1)^\mu} \min_{k\in\iibe{i_\vartheta}{i}}\left(\sum_{j=1}^n
\min\left[\frac{g_{r,k,j}^2}{\varsigma_j},\frac{g_{r,k,j}^2}{v_{r,k}}\right]\right).
\end{align*}
But, using \req{fincrease} in Lemma~\ref{lemma:incrf},
\beqn{fgrowth}
f(x_0)-f(x_{i_\vartheta+1})
= f(x_{r,0})-f(x_{i_\vartheta+1})
\geq -n( \beta_{2,r}+ \half\alpha^2 L)\sum_{k=0}^{i_\vartheta} a(k)^2.
\eeqn
Combining these two last inequalities with AS.3 then gives that
\begin{align*}
\min_{k\in\iibe{i_\vartheta}{i}}\left(\sum_{j=1}^n
\min\left[\frac{g_{r,k,j}^2}{\varsigma_j},\frac{g_{r,k,j}^2}{v_{r,k}}\right]\right)
& \le \frac{2(i+1)^\mu}{\vartheta(i-i_\vartheta)} \Big(f(x_{i_\vartheta+1})-f(x_{i+1})\Big)\\
&  =  \frac{2(i+1)^\mu}{\vartheta(i-i_\vartheta)} \Big(f(x_0)-f(x_{i+1})+f(x_{i_\vartheta+1})-f(x_0)\Big)\\
& \le \frac{2(i+1)^\mu}{\vartheta(i-i_\vartheta)} \Big( \Gamma_0
       + n( \beta_{2,r}+ \half\alpha^2 L)\sum_{k=0}^{i_\vartheta} a(k)^2\Big)
\end{align*}
from which we obtain that there exists a subsequence $\{i_t\}\subseteq
\{i\}_{i_\vartheta}^\infty$ such that
\beqn{A1}
\sum_{j=1}^n
\min\left[\frac{g_{r,i_t,j}^2}{\varsigma_j},\frac{g_{r,i_t,j}^2}{v_{r,i_t}}\right]
\le \frac{2(i_t+1)^\mu}{\vartheta(i_t-i_\vartheta)} \left[ \Gamma_0
       + n( \beta_{2,r}+ \half\alpha^2 L)\sum_{k=0}^{i_\vartheta} a(k)^2\right].
\eeqn
Now,
\beqn{A0}
\frac{(i_t+1)^\mu}{i_t-i_\vartheta}
< \frac{2^\mu i_t^\mu}{i_t-i_\vartheta}
< \frac{2 i_t^\mu}{i_t-i_\vartheta}
=\frac{2 i_t^\mu i_t}{(i_t-i_\vartheta)i_t}
=\frac{i_t}{i_t-i_\vartheta} \cdot \frac{2}{i_t^{1-\mu}}
\le  \frac{2(i_\vartheta+1)}{i_t^{1-\mu}},
\eeqn
where we used the facts that $\mu <1$ and that $\frac{i_t}{i_t-i_\vartheta} $
is a decreasing function for $i_t \ge i_\vartheta+1$. As a
consequence, we obtain from \req{A1} that
\[
\sum_{j=1}^n
\min\left[\frac{g_{r,i_t,j}^2}{\varsigma_j},\frac{g_{r,i_t,j}^2}{v_{r,i_t}}\right]
\le  \frac{4(i_\vartheta+1)}{\vartheta i_t^{1-\mu}}\left[ \Gamma_0
       + n(\beta_{2,r}+ \half\alpha^2 L)\sum_{k=0}^{i_\vartheta} a(k)^2\right].
\]
If we now define
\[
i_\varsigma \eqdef
\left(\frac{4(i_\vartheta+1)\left[\Gamma_0
  +\half n (\beta_{2,r}+ \half\alpha^2 L)\sum_{k=0}^{i_\vartheta} a(k)^2\right]}
     {\vartheta\varsigma_{\min}}\right)^\frac{1}{1-\mu},
\]
we see that, for all $i_t\ge i_\varsigma$,
\[
\min\left[\frac{g_{r,i_t,j}^2}{\varsigma_j},\frac{g_{r,i_t,j}^2}{v_{r,i_t}}\right]
\le \varsigma_{\min}.
\]
We then apply Lemma~\ref{lemma:divsprop} to deduce from \req{A1} that,
for all $i_t\ge i_\varsigma$,
\[
\sum_{j=1}^n\frac{g_{r,i_t,j}^2}{2\varsigma_j}
\le \frac{2(i_t+1)^\mu}{\vartheta(i_t-i_\vartheta)} \left[ \Gamma_0
       + n( \beta_{2,r}+ \half\alpha^2 L)\sum_{k=0}^{i_\vartheta} a(k)^2\right]
\]
and therefore, because $\varsigma_j\le 1$, that
\[
\|g_{r,i_t}\|^2 \
\le \frac{(i_t+1)^\mu}{i_t-i_\vartheta}\left(\frac{2}{\vartheta}\right)\left[ \Gamma_0
       + n( \beta_{2,r}+ \half\alpha^2 L)\sum_{k=0}^{i_\vartheta} a(k)^2\right]
\]
for all $i_t\ge i_\varsigma$. This then gives \req{rate-divs}, the
last inequality following from \req{A0} and \req{fishr}.
\epr

\section{Proof of Theorem~\ref{theorem:adag}}\label{proof-adag}

Again, we first recall a useful technical lemma.

\llem{gen:series}{\cite[Lemma~3.1]{GratJeraToin22b}
Let $\{a_k\}_{k\ge 0}$ be a non-negative sequence,
$\alpha > 0$, $\xi>0$ and define, for each $k \geq 0$,
$b_k = \sum_{t=0}^k a_t$.  Then 
\beqn{alphaall}
\sum_{k=0}^i  \frac{a_k}{(\xi+b_k)^{\alpha}}
\le \left\{\begin{array}{ll}
\bigfrac{1}{(1-\alpha)} ( (\xi + b_i)^{1-  \alpha} - \xi^{1-  \alpha} )
&\tim{if } \alpha \neq 1,\\*[2ex]
\log\left(\bigfrac{\xi + b_i}{\xi} \right)
&\tim{if } \alpha = 1.
\end{array}\right.
\eeqn
}

\paragraph*{Proof of Theorem~\ref{theorem:adag}}
We first note that \req{ws-adag}  implies that 
\beqn{eq3.1}
\varsigma^\mu \leq \max_{j\in\ii{n}}w_{r,i,j} \leq 
\left(\varsigma + \sum_{k=0}^i\|g_{r,k}\|^2\right)^\mu
\eeqn
for all $i\geq0$. We also deduce from \req{fishr} and
\req{fdecrease} in Lemma~\ref{lemma:decrf} with $\ell = r$ that, for $i\geq 0$,
\beqn{summing-up1}
f(x_{i+1})
\le f(x_0) -\beta_{1,r} \sum_{k=0}^i\bigfrac{\|g_{r,k}\|^2}{\max_{j\in\ii{n}}w_{r,k,j}}
+  \left(\beta_{2,r} + \half L\right) 
\bigsum_{j=1}^n\bigsum_{k=0}^i \frac{g_{r,k,j}^2}{w_{r,k,j}^2}.
\eeqn
\paragraph*{Case (i). }
Suppose first that $\mu \in (0,\half)$.
For each $j\in\ii{n}$, we then apply 
\req{alphaall} in Lemma~\ref{gen:series}
with $a_k = g_{r,k,j}^2$, $\xi=\varsigma$ and $\alpha= 2 \mu < 1$, and obtain
from \req{Delta-def} and \req{ws-adag} that,
\beqn{harsh-bound}
\bigsum_{k=0}^i \frac{g_{r,k,j}^2}{w_{r,k,j}^2}
\le \bigfrac{1}{1-2 \mu} \left[\left(\varsigma + \bigsum_{k=0}^i g_{r,k,j}^2\right)^{1-2\mu}
  - \varsigma^{1-2 \mu} \right]
\le  \bigfrac{1}{1-2 \mu} \left(\varsigma + \bigsum_{k=0}^i
g_{r,k,j}^2\right)^{1-2\mu}.
\eeqn
Now substituting this bound in \req{summing-up1} and using AS.3 gives that
\beqn{eq01}
\beta_{1,r}\sum_{k=0}^i \bigfrac{\|g_{r,k}\|^2}{\max_{j\in\ii{n}}w_{r,k,j}}
\le \Gamma_0 + \frac{n(\beta_{2,r}+ \half L)}{1-2\mu}\left(\varsigma + \sum_{k=0}^i \|g_{r,k}\|^2\right)^{1-2\mu}.
\eeqn
Suppose now that
\beqn{hypbad}
\sum_{k=0}^i \|g_{r,k}\|^2
\geq \max\left\{
\varsigma,\frac{1}{2}\left[\frac{(1-2\mu)\Gamma_0}
                           {n(\beta_{2,r}+\half L)}\right]^{\sfrac{1}{1-2\mu}} \right\},
\eeqn
implying
\[
\varsigma + \sum_{k=0}^i \|g_{r,k}\|^2 \leq 2\sum_{k=0}^i \|g_{r,k}\|^2
\tim{and}
\Gamma_0 \leq \frac{n(\beta_{2,r}+\half L)}{1-2\mu}\left(2\sum_{k=0}^i \|g_{r,k}\|^2\right)^{1-2\mu}.
\]
Then, using \req{eq01} and \req{eq3.1},
\[
\frac{\beta_{1,r}}{2^\mu\left[ \sum_{k=0}^i \|g_{r,k}\|^2\right]^\mu}\sum_{k=0}^i \|g_{r,k}\|^2
\leq \frac{2^{2(1-\mu)}n(\beta_{2,r}+\half L)}{1-2\mu} \left(\sum_{k=0}^i \|g_{r,k}\|^2\right)^{1-2\mu}.
\]
Solving this inequality for $\sum_{k=0}^i \|g_{r,k}\|^2$ and using the
fact that $2^{2(1-\mu)}<4$ gives that
\[
\sum_{k=0}^i \|g_{r,k}\|^2 < 
\left[ \frac{4n(\beta_{2,r}+\half L)}{\beta_{1,r}(1-2\mu)}\right]^{\sfrac{1}{\mu}}
\]
and therefore
\beqn{ab1}
\average_{k\in \iiz{i}} \|g_{r,k}\|^2
< \left[ \frac{4n(\beta_{2,r}+\half L)}{\beta_{1,r}(1-2\mu)}\right]^{\sfrac{1}{\mu}}\cdot\frac{1}{{i+1}}.
\eeqn
Alternatively, if \req{hypbad} fails, then
\beqn{ab2}
\average_{k\in\iiz{i}} \|g_{r,k}\|^2
< \max\left\{
\varsigma,\frac{1}{2}\left[\frac{(1-2\mu)\Gamma_0}
                           {n(\beta_{2,r}+\half L)}\right]^{\sfrac{1}{1-2\mu}} \right\}\cdot\frac{1}{i+1}.
\eeqn
Combining \req{ab1}, \req{ab2} and \req{fishr} gives \req{rate-adag} for $0<\mu<\half$.
\paragraph*{Case (ii). }
Let us now consider the case where $\mu = \half$.
For each $j\in\ii{n}$, we apply  
\req{alphaall} in Lemma~\ref{gen:series}
with $a_k = g_{r,k,j}^2$, $\xi = \varsigma$ and $\alpha= 2 \mu = 1$ and obtain that,
\[
\bigsum_{j=1}^n\bigsum_{k=0}^i \frac{g_{r,k,j}^2}{w_{r,k,j}^2}
\le \bigsum_{j=1}^n  \log\left(\frac{1}{\varsigma}\left(\varsigma+ \bigsum_{k=0}^i g_{r,k,j}^2\right)\right)
\le n\log\left(1 + \frac{1}{\varsigma} \bigsum_{k=0}^i \| g_{r,k}\|^2\right)
\]
and substituting this bound in \eqref{summing-up1} then gives that
\[
\beta_{1,r} \sum_{k=0}^i \bigfrac{\|g_{r,k}\|^2}{\max_{j\in\ii{n}}w_{r,k,j}}
\le \Gamma_0 + \half n(\beta_{2,r} + \half L)
\log\left( 1 +  \frac{1}{\varsigma}\bigsum_{k=0}^i \| g_{r,k}\|^2\right).
\]
Suppose now that 
\beqn{logmuhalfhyp}
\sum_{k=0}^i \|g_{r,k}\|^2
\geq \max\left[ \varsigma, e^{\frac{2 \Gamma_0}{n(\beta_{2,r} + \half L)}}\right],
\eeqn
\noindent
implying that 
\[
\varsigma + \sum_{k=0}^i \|g_{r,k}\|^2 \leq 2 \sum_{k=0}^i \|g_{r,k}\|^2
\tim{and}
\Gamma_0 \leq \half n(\beta_{r,2} + \half L)
\log\left( \frac{2}{\varsigma}\bigsum_{k=0}^i \| g_{r,k}\|^2\right).
\]
\noindent
Using \req{eq3.1} for $\mu=\half$, we obtain then that
\[
\bigfrac{\beta_{1,r}}{\sqrt{2} \,\sqrt{\sum_{k=0}^i \| g_{r,k}\|^2}} \bigsum_{k=0}^i  \|g_{r,k}\|^2
\le  n(\beta_{r,k} +\half L)
\vspace*{-2mm}
\]
and thus that
\beqn{tosolvmuhalf}
 \beta_{r,1} \sqrt{\bigsum_{k=0}^i \| g_{r,k}\|^2}
 \le 2\sqrt{2}\max[3\beta_{r,1}, n(\beta_{r,2} + \half L)]
 \log\left( \sqrt{\frac{2}{\varsigma} \bigsum_{k=0}^i \| g_{r,k}\|^2} \right).
\eeqn
Now define
\beqn{ggu}
\gamma_1 \eqdef \beta_{r,1}\sqrt{\frac{ \varsigma}{2}},
\ms
\gamma_2 \eqdef  2\sqrt{2} \max[\sfrac{3}{2}\beta_{r,1}, n(\beta_{r,2} + \half L)] 
\tim{ and }
u \eqdef \sqrt{\frac{2}{\varsigma} \bigsum_{k=0}^i\|g_{r,k}\|^2}
\eeqn
and observe that $\gamma_2 > 3 \gamma_1$ by construction, because $\varsigma\leq 1$. 
The inequality \req{tosolvmuhalf} can then be rewritten as
\beqn{tosolvemuhalf}
\gamma_1 u \le \gamma_2 \log(u).
\eeqn
Let us denote by $\psi(u) \eqdef \gamma_1 u - \gamma_2 \log(u)$. Since
$\gamma_2 > 3 \gamma_1$, the equation $\psi(u)=0$ admits two roots $u_1 \leq
u_2$ and \req{tosolvemuhalf} holds for $u\in[u_1,u_2]$.
The definition of $u_2$ then gives that
\[
\log(u_2)- \frac{\gamma_1}{\gamma_2}u_2 = 0.
\]
Setting $z = -\frac{\gamma_1}{\gamma_2}u_2$, we obtain that
\[
z e^z = -\frac{\gamma_1}{\gamma_2}.
\]
Thus $z = W_{-1}(-\frac{\gamma_1}{\gamma_2})<0$, where $W_{-1}$ is the second
branch of the Lambert function defined over $[-\frac{1}{e}, 0)$.
As $-\frac{\gamma_1}{\gamma_2} \geq -\frac{1}{3} $, $z$ is well defined and thus
\[
u_2
= -\frac{\gamma_2}{\gamma_1}\,z
= -\frac{\gamma_2}{\gamma_1}\,W_{-1}\left(-\frac{\gamma_1}{\gamma_2}\right)>0
= -\psi\,W_{-1}\left(-\frac{1}{\psi}\right)>0,
\]
where $\psi = \frac{\gamma_2}{\gamma_1}$ is given by \req{psi-def}.
As a consequence, we deduce from \req{tosolvemuhalf} and \req{ggu} that
\[
\bigsum_{k=0}^i\|g_{r,k}\|^2
= \frac{\varsigma}{2}\,u_2^2
=\frac{\varsigma\psi^2}{2} \,\left|W_{-1}\left(-\frac{1}{\psi}\right)\right|^2
\]
and
\beqn{firstmuhalfineq}
\average_{k\in\iiz{i}} \|g_{r,k}\|^2
\leq \frac{\varsigma\psi^2}{2} \,\left|W_{-1}\left(-\frac{1}{\psi}\right)\right|^2
\cdot\frac{1}{i+1}.
\eeqn
If \eqref{logmuhalfhyp} does not hold, we have that
\beqn{secmuhalfineq}
\average_{k\in\iiz{i}} \|g_{r,k}\|^2
< \max\left[ \varsigma, e^{\frac{2 \Gamma_0}{n(\beta_{2,r} + \half L)}}\right]\cdot \frac{1}{i+1}.
\eeqn
Combining \req{firstmuhalfineq}, \req{secmuhalfineq} and \req{fishr}
gives \req{rate-adag} for $\mu=\half$.
\paragraph*{Case (iii). }
Finally, suppose that $\half < \mu < 1$.
Once more, we apply 
\req{alphaall} in Lemma~\ref{gen:series} for each $j\in\ii{n}$
with $a_k = g_{r,k,j}^2$, $\xi=\varsigma$ and $\alpha= 2 \mu > 1$ and obtain that, for $j\in\ii{n}$,
\beqn{mugeqhalfound}
\bigsum_{k=0}^i \frac{g_{r,k,j}^2}{w_{r,k,j}^2}
\le \frac{1}{1-2 \mu} \left( \Big(\varsigma + \sum_{k=0}^i g_{r,k,j}^2\Big)^{1-2\mu}- \varsigma^{1-2\mu} \right)
\le \frac{\varsigma^{1-2 \mu}}{2\mu -1}.
\eeqn
Substituting the bound \req{mugeqhalfound} in \req{summing-up1} and using
\req{eq3.1} and AS.3 gives that
\[
\beta_{r,1}\sum_{k=0}^i  \frac{1}{(\varsigma +
\sum_{t=0}^i \| g_{r,t}\|^2 )^\mu}\|g_{r,k}\|^2 
\le \Gamma_0 + \frac{n(\beta_{2,r} +\half  L)\varsigma^{1-2\mu}}{2\mu-1}. \\
\]
If we now suppose that
\beqn{musuphalfhyp}
\sum_{k=0}^i \| g_{r,k}\|^2 \geq \varsigma,
\eeqn
then
\beqn{firstmugeqhalf}
\average_{k\in\iiz{i}} \|g_{r,k}\|^2
\le \left[ \frac{2^\mu}{\beta_{r,1}} \left(\Gamma_0
  + \frac{n(\beta_{2,r} + \half L)\varsigma^{1-2\mu}}{2\mu-1}\right) \right]^{\sfrac{1}{1-\mu}}
  \cdot \frac{1}{i+1}.
\eeqn
If \req{musuphalfhyp} does not hold, we derive that
\beqn{secmugeqhalf}
\average_{k\in\iiz{i}} \|g_{r,k}\|^2 \leq \frac{\varsigma}{(i+1)}.
\eeqn
Thus, \req{firstmugeqhalf}, \req{secmugeqhalf} and \req{fishr} finally
imply \req{rate-adag} for $\half<\mu<1$. 
\epr

\end{document}